\documentclass[10pt,reqno]{amsproc}

\usepackage{amssymb}
\usepackage[mathscr]{eucal}
\usepackage{enumitem}

\usepackage[all]{xy}
\ifx\pdfoutput\undefined \xyoption{dvips}\fi

\newtheorem{thm}{Theorem}
\newtheorem{lemma}[thm]{Lemma}
\newtheorem{cor}[thm]{Corollary}

\theoremstyle{definition}
\newtheorem{defn}[thm]{Definition}
\newtheorem{notation}[thm]{Notation}
\newtheorem{obs}[thm]{Observation}
\newtheorem{recall}[thm]{Recall}
\newtheorem{eg}[thm]{Example}

\def\parhead[#1]{\vspace{1ex}%
\noindent{\bf\boldmath #1.}}

\newcommand{\mlaux}[3]{\setbox0=\hbox{$\mathsurround=0pt #2{#3}$}%
  \dimen0=\dp0\advance\dimen0 by \ht0\lower#1\dimen0\box0}
\newcommand{\mlower}[2]{\mathpalette{\mlaux{#1}}{#2}}

\newcommand{\makekdp}[2]{\setbox0=\hbox{$\mathsurround=0pt #1{#2}$} \dp0=0pt \box0}
\newcommand{\killdepth}{\relax\mathpalette\makekdp}
\newcommand{\category}[1]{\underline{\killdepth{\text{#1}}}}

\newcommand{\makellapm}[2]{\hbox to 0pt{\hss$\mathsurround=0pt #1{#2}$}}

\newcommand{\makerlapm}[2]{\hbox to 0pt{$\mathsurround=0pt #1{#2}$\hss}}

\newcommand{\icat}{\mathbb}
\newcommand{\lcat}{\mathcal}

\newcommand{\defeq}{\mathrel{\mlower{0.15}{\stackrel{\text{def}}=}}}

\newcommand{\ifundef}[1]{\expandafter\ifx\csname#1\endcsname\relax}

\ifundef{Top}
  \DeclareSymbolFont{dsymbolsC}{U}{txsyc}{m}{n}
  \SetSymbolFont{dsymbolsC}{bold}{U}{txsyc}{bx}{n}
  \DeclareFontSubstitution{U}{txsyc}{m}{n}

  \def\re@DeclareMathSymbol#1#2#3#4{%
      \let#1=\undefined
      \DeclareMathSymbol{#1}{#2}{#3}{#4}}

  \re@DeclareMathSymbol{\Top}{\mathord}{dsymbolsC}{120}
  \re@DeclareMathSymbol{\Bot}{\mathord}{dsymbolsC}{121}
\fi

\newcommand{\partinj}{\Bot}
\newcommand{\partproj}{\Top}

\let\newmop=\DeclareMathOperator

\newmop{\im}{im}
\newmop{\coim}{coim}

\newmop{\cell}{cell}
\newmop{\cof}{cof}
\newmop{\fib}{fib}

\newmop{\triv}{th}
\newmop{\sst}{sp}
\newmop{\sk}{sk}
\newmop{\cosk}{ck}
\newmop{\dec}{dec}
\newmop{\lax}{lax}

\newmop{\cls}{cls}

\newmop{\depth}{dp}
\newmop{\id}{id}

\newmop{\colim}{colim}

\newmop{\obj}{obj}
\newmop{\arr}{arr}
\newmop{\sq}{sq}

\newmop{\dis}{dis}

\newmop{\cod}{cod}
\newmop{\dom}{dom}

\newmop{\rev}{rev}

\newcommand{\dual}{^\circ}
\newcommand{\oth}{^{\mathord{\text{th}}}}

\newcommand{\op}{^{\mathord{\text{\rm op}}}}

\newcommand{\boundary}{\partial}

\newcommand{\pretens}{\mathbin{\boxtimes}}
\newcommand{\spretens}{\mathbin{\boxdot}}

\newcommand{\gray}{\mathbin{\circledast}}
\newcommand{\laxgray}{\mathbin{\otimes}}

\newcommand{\ctimes}{\times_c}

\newcommand{\cpretens}{\pretens_c}

\newcommand{\cgray}{\gray_c}

\newcommand{\sgray}{_{\gray}}
\newcommand{\slgray}{_{\laxgray}}

\newcommand{\Slgray}{{\Strat\slgray}}
\newcommand{\hclgray}{\hcpath\slgray}

\newcommand{\pspath}{\bar{\hcpath}}
\newcommand{\hcpath}{\icat{S}}
\newcommand{\hcpathg}{\hat{\hcpath}}
\newcommand{\hchorn}{\mathbb{H}}

\newcommand{\indec}{}
\def\indec<#1,#2>{\langle#1,#2\rangle}

\newcommand{\sint}{}
\newcommand{\hint}{}
\def\sint(#1,#2){{(#1,#2)}}
\def\hint(#1,#2]{{(#1,#2]}}

\let\newfunctor=\newmop

\newfunctor{\Dec}{Dec}

\newcommand{\inerve}{\nerv_\infc}
\newcommand{\inladj}{\ladj_\infc}

\newfunctor{\pathecat}{P}
\newfunctor{\recons}{R}

\newfunctor{\forget}{U}
\newfunctor{\refl}{L}
\newfunctor{\dgm}{D}
\newfunctor{\ladj}{F}
\newfunctor{\gfunc}{G}
\newfunctor{\incl}{I}
\newfunctor{\nerv}{N}

\newfunctor{\efunc}{E}
\newfunctor{\kfunc}{K}
\newfunctor{\lfunc}{L}

\newcommand{\yoneda}[1]{\ulcorner{#1}\urcorner}

\newcommand{\satom}[1]{\lbrack\!\lbrack{#1}\rbrack\!\rbrack}
\newcommand{\statom}{\satom}

\newdir{ >}{{}*!/-10pt/@{>}}
\newdir{u(}{{}*!/-5pt/@^{(}}
\newdir{d(}{{}*!/-5pt/@_{(}}
\newdir{|>}{%
  !/4.5pt/@{|}*:(1,-.2)@^{>}*:(1,+.2)@_{>}*!/-3pt/@{ }}

\def\arrow#1:#2->#3.{{#1\colon #2}\xy (0,0)*+{} %
  \ar (7,0)*+{}\endxy{#3}}

\def\epi#1:#2->#3.{{#1\colon #2}\xy (0,0)*+{} %
  \ar@{->>} (7,0)*+{}\endxy{#3}}

\def\overepi#1:#2->#3.{{#2}\xy (0,0)*+{} %
  \ar@{->>}^-{\smash{#1}} (7,0)*+{}\endxy{#3}}

\def\cover#1:#2->#3.{{#1\colon #2}\xy (0,0)*+{} %
  \ar@{-|>} (7,0)*+{}\endxy{#3}}

\def\overinc#1:#2->#3.{{#2}\xy (0,0)*+{} %
  \ar@{u(->}^-{\smash{#1}} (9,0)*+{}\endxy{#3}}

\def\overarr#1:#2->#3.{{#2}\xy (0,0)*+{} %
  \ar^-{\smash{#1}} (7,0)*+{}\endxy{#3}}

\def\spanarr#1:#2->#3.{{#1\colon #2}\xy (0,0)*+{} %
  \ar|-{\object@{|}} (9,0)*+{}\endxy{#3}}

\def\inc#1:#2->#3.{{#1\colon #2}\xy (0,0)*+{} %
  \ar@{u(->} (9,0)*+{}\endxy{#3}}

\def\nattrans#1:#2->#3.{{#1\colon #2}\xy (0,0)*+{} %
  \ar^-{.} (7,0)*+{}\endxy{#3}}

\def\iso#1:#2->#3.{{#1\colon #2}\xy (0,0)*+{} %
  \ar^-{\smash{\simeq}} (7,0)*+{}\endxy{#3}}

\def\equiv#1:#2->#3.{{#1\colon #2}\xy (0,0)*+{} %
  \ar^-{\smash{\sim}} (7,0)*+{}\endxy{#3}}

\def\anoniso#1->#2.{{#1}\xy (0,0)*+{} %
  \ar^-{\simeq} (7,0)*+{}\endxy{#2}}

\def\twocell#1:#2->#3.{{#1\colon #2}\xy (0,0)*+{} %
  \ar@{=>} (7,0)*+{}\endxy{#3}}

\def\isotwocell#1:#2->#3.{{#1\colon #2}\xy (0,0)*+{} %
  \ar@{=>}^-{\sim} (7,0)*+{}\endxy{#3}}

\newcommand{\Wcs}{\category{Wcs}}
\newcommand{\Gray}{\category{Gray}}
\newcommand{\Grph}{\category{Grph}}

\newcommand{\ECat}[1]{{#1}\text{-}\Cat}

\newcommand{\Simp}{\category{Simp}}
\newcommand{\Set}{\category{Set}}

\newcommand{\Cat}{\category{Cat}}

\newcommand{\Strat}{\category{Strat}}

\newcommand{\infc}{\ensuremath\omega}
\newcommand{\InfCat}{\ensuremath{\ECat\infc}}

\newcommand{\aDelta}{\Delta_{\mathord{+}}}

\newcommand{\face}{\delta}
\newcommand{\vertex}{\varepsilon}
\newcommand{\degen}{\sigma}

\newcommand{\cube}[1]{{[\![#1]\!]}}

\newcommand{\pocorner}{\hbox to 10pt{{\vrule height10pt depth0pt width0.5pt}%
    \vbox to 10pt{{\hrule height0.5pt width9.5pt depth0pt}\vfill}}}
\newcommand{\poexcursion}{\save[]-<12pt,-12pt>*{\pocorner}\restore}
\newcommand{\pbcorner}{\vbox to 0pt{\kern 5pt\hbox to 0pt{\kern 5pt%
      \vbox{{\hrule height0.5pt width9.5pt depth0pt}}%
      {\vrule height10pt depth0pt width0.5pt}\hss}\vss}}

\def\funcdisplay #1:#2->#3. {{%
    \xymatrix@R=0em@C=14em{%
      {#2}\ar[r]^-{#1} & {#3}}}}

\def\adjdisplay#1->#2,#3-|#4.{{%
    \xymatrix@R=0em@C=7em{%
      {#1} \ar@/_0.65pc/[rr]_-{#4} & {\bot} &
      {#2}\ar@/_0.65pc/[ll]_-{#3}}}}

\def\adjinline#1->#2,#3-|#4.{\arrow {#3}\dashv{#4}:#1->#2.}

\newcommand{\pent}[2]{\xybox{\xymatrix@!=#2{
    & {\c}\ar[rr]^{\cd} && {\d}\ar[rdd]^{\de} & \\
    &&&& \\
    {\b}\ar[uur]^{\bc} &&&& {\e} \\
    && {\a}\ar[ull]^{\ab}\ar[urr]_{\ae} && {\ifcase #1
      \ar"4,3";"1,2"|*+{\labelstyle \ac}="one"
      \ar"4,3";"1,4"|*+{\labelstyle \ad}="two"
      \ar@{}"3,1";"one"|(0.55){\objectstyle \abc}
      \ar@{}@<1.5em>"one";"two"|{\objectstyle \acd}
      \ar@{}"two";"3,5"|(0.45){\objectstyle \ade} \or
      \ar"3,1";"1,4"|*+{\labelstyle \bd}="one"
      \ar"3,1";"3,5"|*+{\labelstyle \be}="two"
      \ar@{}"1,2";"one"|(0.55){\objectstyle \bcd}
      \ar@{}@<1.5em>"one";"two"|{\objectstyle \bde}
      \ar@{}"two";"4,3"|(0.45){\objectstyle \abe} \or
      \ar"4,3";"1,2"|*+{\labelstyle \ac}="two"
      \ar"1,2";"3,5"|*+{\labelstyle \ce}="one"
      \ar@{}"1,4";"one"|(0.55){\objectstyle \cde}
      \ar@{}@<1.5em>"one";"two"|{\objectstyle \ace}
      \ar@{}"two";"3,1"|(0.45){\objectstyle \abc} \or
      \ar"3,1";"1,4"|*+{\labelstyle \bd}="two"
      \ar"4,3";"1,4"|*+{\labelstyle \ad}="one"
      \ar@{}"3,5";"one"|(0.55){\objectstyle \ade}
      \ar@{}@<1.5em>"one";"two"|{\objectstyle \abd}
      \ar@{}"two";"1,2"|(0.45){\objectstyle \bcd} \or
      \ar"3,1";"3,5"|*+{\labelstyle \be}="one"
      \ar"1,2";"3,5"|*+{\labelstyle \ce}="two"
      \ar@{}"4,3";"one"|(0.55){\objectstyle \abe}
      \ar@{}@<1.5em>"one";"two"|{\objectstyle \bce}
      \ar@{}"two";"1,4"|(0.45){\objectstyle \cde} \else\fi} }}}
\newcommand{\pentofpent}[1]{
  \def\baselen{#1}
  \begin{xy}
    0;<\baselen,0mm>:
    *{\xybox{
        \POS(0,-4)*[o]{\pent 0{\baselen}}="zero"
        \POS(16,40)*[o]{\pent 3{\baselen}}="three"
        \POS(72,40)*[o]{\pent 1{\baselen}}="one"
        \POS(88,-4)*[o]{\pent 4{\baselen}}="four"
        \POS(44,-36)*[o]{\pent 2{\baselen}}="two"
        \ar@<1ex>"zero";"three"^-{\objectstyle\abcd}
        \ar@<1ex>"three";"one"^-{\objectstyle\abde}
        \ar@<1ex>"one";"four"^-{\objectstyle\bcde}
        \ar@<-1ex>"zero";"two"_-{\objectstyle\acde}
        \ar@<-1ex>"two";"four"_-{\objectstyle\abce}
        \ar(44,-5);(44,+15)^{\objectstyle\abcde}
     }}
  \end{xy}
}

\input{parti_labels.aux}

\title[Weak Complicial Sets II]
      { Weak Complicial Sets \\ 
        A Simplicial Weak \infc-Category Theory \\
        Part II: Nerves of Complicial Gray-Categories
}

\author[Verity]{Dominic Verity}
\dedicatory{To Ross Street on the occasion of his $60\oth$ birthday.}
\address{
  Centre of Australian Category Theory \\
  Macquarie University \\
  NSW 2109 \\
  Australia
}
\email{domv@ics.mq.edu.au}

\thanks{The author would like to thank Fitzwilliam College Cambridge
  whose visiting fellowship programme supported the completion of this
  work.}

\date{March 1, 2006}


\begin{document}

\begin{abstract}
  This paper continues the development of a simplicial theory of weak
  \infc-categories, by studying categories which are enriched in {\em
    weak complicial sets}. These {\em complicial Gray-categories\/}
  generalise both the Kan complex enriched categories of homotopy
  theory and the 3-categorical Gray-categories of weak 3-category
  theory. We derive a simplicial nerve construction, which is closely
  related to Cordier and Porter's {\em homotopy coherent nerve}, and
  show that this faithfully represents complicial Gray-categories as
  weak complicial sets. The category of weak complicial sets may
  itself be canonically enriched to a complicial Gray-category whose
  homsets are higher generalisations of the bicategory of
  homomorphisms, strong transformations and modifications. By applying
  our nerve construction to this structure, we demonstrate that the
  totality of all (small) weak complicial sets and their structural
  morphisms at higher dimensions form a richly structured (large) weak
  complicial set.
\end{abstract}

\maketitle
\tableofcontents


\section{Introduction}
\label{intro.sec.a}

This paper is the second in a series of works on the simplicial weak
\infc-category theory of {\em weak complicial sets}. In the first of
these~\cite{Verity:2005:WeakComp} we studied the fundamental homotopy
theory of such structures by constructing a canonical model structure
on the category of {\em stratified simplicial sets\/} whose fibrant
objects are these weak complicial sets. On the way we developed a
range of technical tools for studying such objects, including an
equivalence based characterisation of outer horn filling, Gray tensor
products and associated complicial closures, join and d{\'e}calage
constructions and so forth. We also showed that weak complicial sets
usefully subsume and generalise a range of common weakened and higher
categorical notions, including Kan complexes, Joyal's quasi-categories
and strict \infc-categories.

Having proposed a candidate simplicial weak \infc-category theory we
might now wish to establish its credentials from a number of
perspectives. From {\em below}, we should be interested in
establishing the bounds of its expressiveness by demonstrating that it
may be used to faithfully represent a range of common weak
\infc-categorical notions and their ``strictifications''. At the very
least, we should be prepared to demonstrate how the well understood
low dimensional cases may be represented in this way. We might also
wish to undertake a similar program with regard to higher strict
structures such as the {\em classical Gray-categories}, that is to say
categories enriched in strict \infc-categories with respect to their
Gray tensor product. From {\em above}, we must also have regard to
ensuring that the devices we define carry enough structure to allow us
to prove strong structural theorems. Category theorists often code
this intuition into {\em coherence theorems\/} demonstrating that our
weakened structures may be {\em rectified\/} and replaced by
equivalent strict structures of a certain class.

Of course the completion of such an investigation would not, in itself,
establish the full utility of our theory. For example, this purely
structural study might completely circumvent the development of the
actual category theory of these structures. However, if we are to use
these objects in practice they must come equipped with all of the
basic constructions of ($n$-)category theory, such as the
(generalised) span construction. We should also, at the very least, be
able to establish homotopical versions of all of the standard results,
constructions and structures of classical category theory, such as
discrete and categorical fibrations, Yoneda's lemma, adjunctions,
limits and colimits and so forth. This latter project is the subject
of another companion paper~\cite{Verity:2006:IntQCat}, wherein we
represent weak complicial sets as certain kinds of {\em complicially
  enriched quasi-categories}. This provides us with a natural context
in which to generalise traditional category theory to a kind of
homotopy coherent quasi-category theory within the Quillen model
category of weak complicial sets itself. 

Our far more modest goal here is simply to take the first steps toward
implementing the coherence and representation program of the paragraph
before last. The completion of that program will have to wait for
later papers in this series. Herein we generalise the Gray-categories
of weak 3-category theory to {\em complicial Gray-categories}, which
are simply categories enriched in weak complicial sets with respect to
the Gray tensor product of stratified sets. Examples of such
structures include classical Gray-categories and the category of weak
complicial sets itself. Then we extend the homotopy coherent nerve
construction of Cordier and Porter~\cite{Cordier:1986:HtyCoh} to
provide a nerve functor which faithfully represents complicial
Gray-categories as weak complicial sets. Finally we apply this to
demonstrate that the totality of all (small) weak complicial sets and
their structural morphisms may be gathered together to form a richly
structured (large) weak complicial set. 

In later work, our eventual goal will be to prove a coherence theorem,
which we might state informally as:
\begin{quotation}
  \em Every weak complicial set satisfying certain mild {\em
    equivalence repletion\/} conditions on its thin simplices is
  homotopy equivalent to the nerve of some complicial Gray-category
  whose homsets also satisfy that condition.
\end{quotation}
It is these structures amongst the weak complicial sets which we might
truly call our simplicial weak \infc-categories.


\section{Background Definitions and Notational Conventions}

This paper should be read in conjunction
with~\cite{Verity:2005:WeakComp}. In particular, we will rely upon
that work to describe the notational conventions used here and to
provide a foundational introduction to the homotopy theory of weak
complicial sets.  However, as a warm up we recall a few basic
definitions and pivotal results here:

\begin{recall}[simplicial and stratified sets]
  As usual we let $\aDelta$ denote the (skeletal) category of finite
  ordinals and order preserving maps between them and use the notation
  $\Delta$ to denote its full subcategory of non-zero
  ordinals. Following tradition we let $[n]$ denote the ordinal $n+1$
  as an object of $\aDelta$ and refer to arrows of $\aDelta$ as {\em
    simplicial operators\/} which are denoted using Greek letters
  $\alpha,\beta,...$. Then the category $\Simp$ of {\em simplicial
    sets\/} and {\em simplicial maps\/} between them is simply the
  functor category $[\Delta\op,\Set]$, where $\Set$ denotes the
  (large) category of all (small) sets and functions between them.  We
  shall assume that the reader is familiar the basic homotopy theory
  of simplicial sets as expounded in~\cite{GabrielZisman:1967:CFHT}
  or~\cite{May:1967:Simp}. We shall make repeated use of the classical
  nerve construction $\arrow N:\Cat->\Simp.$ which associates with
  each (small) category $\icat{C}$ a simplicial set $N(\icat{C})$
  whose $n$-simplices correspond to composable paths of arrow of
  $\icat{C}$ of length $n$.
  
  A {\em stratification\/} on a simplicial set $X$ is a subset $tX$ of
  its simplices satisfying the conditions that
  \begin{itemize}
  \item no $0$-simplex of $X$ is in $tX$, and
  \item all of the degenerate simplices of $X$ are in $tX$.
  \end{itemize}
  A {\em stratified set\/} is a pair $(X,tX)$ consisting of a
  simplicial set $X$ and a chosen stratification $tX$ the elements of
  which we call {\em thin\/} simplices. In practice, we will elect to
  notionally confuse stratified sets with their underlying simplicial
  sets $X,Y,Z,...$ and uniformly adopt the notation $tX,tY,tZ,...$ for
  corresponding sets of thin simplices.  A {\em stratified map\/}
  $\arrow f:X->Y.$ is simply a simplicial map of underlying simplicial
  sets which {\em preserves thinness\/} in the sense that for all
  $x\in tX$ we have $f(x)\in tY$. Identities and composites of
  stratified maps are clearly stratified maps, from which it follows
  that we have a category $\Strat$ of stratified sets and maps. 

  We refer the reader to~\cite{Verity:2006:Complicial}
  and~\cite{Verity:2005:WeakComp} for a full development of the theory
  of stratified sets. In particular, we should draw attention to
  definition~\ref{reg.ent.def} of the latter which defines what it
  means to be a {\em stratified subset\/} and identifies two important
  classes of such subsets, the {\em regular\/} and {\em entire\/}
  ones. 
\end{recall}

\begin{recall}[elementary operators]
  We will use the following standard notation and nomenclature
  throughout:
  \begin{itemize}
  \item The injective maps in $\aDelta$ are referred to as {\em face
      operators}. For each $j\in[n]$ we use the
    $\arrow\face^n_j:[n-1]->[n].$ to denote the {\em elementary\/}
    face operator distinguished by the fact that its image does
    not contain the integer $j$.
  \item The surjective maps in $\aDelta$ are referred to as {\em
      degeneracy operators}. For each $j\in[n]$ we use
    $\arrow\degen^n_j:[n+1]->[n].$ to denote the {\em elementary\/}
    degeneracy operator determined by the property that two integers
    in its domain map to the integer $j$ in its codomain.
  \item For each $i\in[n]$ the operator $\arrow\vertex^n_i:[0]->[n].$
    given by $\vertex^n_i(0)=i$ is called the {\em $i\oth$
      vertex operator of $[n]$}. 
  \item We also use the notations $\arrow\eta^n:[n]->[0].$ and
    $\arrow\iota^n:[-1]->[n].$ to denote the unique such simplicial
    operators.
  \end{itemize}
  Unless doing so would introduce an ambiguity, we will tend to reduce
  notational clutter by dropping the superscripts of these elementary
  operators.
\end{recall}

\begin{notation}\label{simp.of.1.a}
  We also use the following notation for the simplices of the
  standard simplex $\Delta[1]$:
  \begin{itemize}
  \item $\arrow 0^r:[r]->[1].$ is the operator which maps each
    $i\in[r]$ to $0\in[1]$.
  \item $\arrow 1^r:[r]->[1].$ is the operator which maps each
    $i\in[r]$ to $1\in[1]$.
  \item $\arrow\rho^r_i:[r]->[1].$ ($1\leq i\leq r$) is the operator
    defined by
    \begin{equation*}
      \rho^r_i(j) =
      \begin{cases}
        0 & \text{if $j< i$,} \\
        1 & \text{if $j\geq i$.}
      \end{cases}
    \end{equation*}
  \end{itemize}
  Later on it will become convenient to index the $r$-simplices of
  $\Delta[1]$ using the {\em doubly pointed\/} set
  $\cube{r}\defeq\{-,+,1,2,...,r\}$, by letting $\rho^r_{-}=0^r$,
  $\rho^r_{+}=1^r$ and defining $\rho^r_i$ as above for an arbitrary
  integer (non-point) in $\cube{r}$.
\end{notation}

\begin{recall}[anodyne extensions and weak complicial sets]
  We refer the reader to notations~\ref{stand.simp.def}
  and~\ref{comp.simp.def} of~\cite{Verity:2005:WeakComp} for the
  definitions of the following stratified sets:
  \begin{itemize}
  \item the {\em standard $n$-simplex\/} $\Delta[n]$, its {\em
      boundary\/} $\boundary\Delta[n]$ and the {\em standard thin
      $n$-simplex\/} $\Delta[n]_t$.
  \item the {\em $k$-complicial $n$-simplex \/} $\Delta^k[n]$ and its
    variants $\Delta^k[n]'$ and $\Delta^k[n]''$, and
  \item the {\em $(n-1)$-dimensional $k$-complicial horn\/}
    $\Lambda^k[n]$.
  \end{itemize}
  The set of {\em elementary anodyne extensions\/} in $\Strat$
  consists of two families of subset inclusions:
  \begin{itemize}
  \item $\overinc\subseteq_r:\Lambda^k[n]->\Delta^k[n].$ for
    $n=1,2,...$ and $k\in[n]$, these are called {\em complicial horn
      extensions}, and
  \item $\overinc\subseteq_e:\Delta^k[n]'->\Delta^k[n]''.$ for
    $n=2,3,...$ and $k\in[n]$, these are called {\em complicial
      thinness extensions}.
  \end{itemize}
  We classify these elementary anodyne extensions into two
  sub-classes, the {\em inner\/} ones for which the index $k$
  satisfies $0<k<n$ and the remaining {\em left\/} and {\em right
    outer\/} ones for which $k=0$ or $k=n$ respectively.  We say that
  a stratified inclusion (monomorphism) $\inc e:U->V.\in\Strat$ is an
  {\em (inner) anodyne extension\/} if it is in the cellular
  completion (that is the completion under pushouts and transfinite
  composition) of the set of elementary (inner) anodyne extensions.  

  A stratified set $A$ is said to be a
  \begin{itemize}
  \item {\em weak inner complicial set\/} if it has the right lifting
    property (RLP) with respect to all inner elementary anodyne
    extensions.
  \item {\em weak complicial set\/} if it has the RLP with respect to
    all elementary anodyne extensions.
  \end{itemize}
  Informally we might simply say that a weak complicial set has {\em
    fillers\/} for all complicial horns. We also say that a stratified
  map $\arrow p:X->Y.$ is a {\em (inner) complicial fibration\/} iff
  it has the RLP with respect to elementary (inner) anodyne extensions.
\end{recall}

\begin{recall}[Gray tensor products]
  The category $\Strat$ of stratified sets supports three related
  tensor products:
  \begin{itemize}[label={}, fullwidth, leftmargin=1em, 
    labelsep=0em, itemsep=1ex]
  \item The {\bf Gray tensor product}
    $\arrow\gray:\Strat\times\Strat->\Strat.$ which is simply defined
    to be the cartesian product construction in $\Strat$. Its closure
    is denoted $\hom(X,Y)$ and is often referred to as the {\em
      stratified set of strong transformations\/} (cf.\
    definition~\ref{gray.tensor.def} of~\cite{Verity:2005:WeakComp}).
  \item The {\bf lax Gray tensor product}
    $\arrow\laxgray:\Strat\times\Strat->\Strat.$ which acts like the
    cartesian product on the underlying simplicial sets of stratified
    sets $X$ and $Y$ but which enjoys a stratification which is a
    subset of that of $X\gray Y$
    (cf.\ definition~\ref{laxgray.tensor.def}
    of loc.\ cit.). This tensor makes $\Strat$ into a
    monoidal category but it is neither symmetric nor does it have a
    corresponding closure.
  \item The {\bf lax Gray pre-tensor product}
    $\arrow\pretens:\Strat\times\Strat->\Strat.$ which again acts like
    the cartesian product on underlying simplicial sets and enjoys a
    stratification which is a subset of that of the lax Gray tensor
    $X\laxgray Y$ (cf.\ definition~\ref{pretens.defn} of loc.\
    cit.). Unlike $\laxgray$, this pre-tensor $\pretens$ has left and
    right closures $lax_l(X,Y)$ and $\lax_r(X,Y)$, called the {\em
      stratified sets of left and right lax transformations}, but it
    is not associative so fails to make $\Strat$ into a monoidal
    category. However, the tensors $\laxgray$ and $\pretens$ are
    indistinguishable from the point of view of weak (inner)
    complicial sets, in the sense that each inclusion
    $\overinc\subseteq_e:X\pretens Y->X\laxgray Y.$ is an inner
    anodyne extension (cf. observation~\ref{pretens.gen.obs} of loc.\
    cit.).
  \end{itemize}
  These tensors generalise the Gray tensor products of $2$-category
  theory to the current context. In particular, it is shown
  in~\cite{Verity:2006:Complicial} that the lax Gray tensor $\laxgray$
  coincides with the Gray tensor product of \infc-categories of
  Steiner~\cite{Steiner:1991:Tensor} or Crans~\cite{Crans:1995:PhD}
  under Street's \infc-categorical nerve
  construction~\cite{Street:1987:Oriental}.
\end{recall}

\begin{recall}[the complicial model structure]
  The category $\Strat$ supports a Quillen model structure, called the
  {\em complicial model structure}, under which:
  \begin{itemize}
  \item A stratified map $\arrow w:X->Y.$ is a weak equivalence if and
    only if the stratified map $\arrow\hom(w,A):\hom(Y,A)->\hom(X,A).$
    is a {\em homotopy equivalence\/}, in the sense of
    definition~\ref{hty.equiv.def} of~\cite{Verity:2005:WeakComp}, for
    each weak complicial set $A$.
  \item The cofibrations are precisely the inclusions of stratified
    sets, thus ensuring that all stratified sets are cofibrant.
  \item The fibrant objects are the weak complicial sets.
  \item An arrow $\arrow p:A->B.$ between weak complicial sets is a
    fibration, called a {\em completely complicial fibration}, if and
    only if it is a complicial fibration
    (cf. lemma~\ref{cfib+wcs=>ccfib} of loc.\ cit.).  All
    completely complicial fibrations are actually complicial
    fibrations but the reverse implication does not necessarily hold
    in general.
  \item An inclusion $\inc e:U->V.$ is a trivial cofibration, called a
    {\em complicial cofibration}, if and only if it has the left
    lifting property (LLP) with respect to all complicial fibrations
    between weak complicial sets (cf.\ corollary~\ref{ccof.char.b} of
    loc.\ cit.). In particular, all anodyne extensions are complicial
    cofibrations. 
  \end{itemize}
  Details of the construction of this model structure may be found in
  section~\ref{quillen.model.sec} of loc.\ cit.
\end{recall}



\section{Complicial Gray-Categories}

We start with a definition, which we claim is the appropriate
generalisation of the tricategorical notion of {\em Gray-category\/}
to our current context.

\begin{defn}[$\Strat\sgray$-categories and
  $\Slgray$-categories]
  A {\em $\Strat\sgray$-category\/} is a category which is enriched
  over the category of stratified set $\Strat$ with respect to the
  Gray tensor product $\gray$. Correspondingly, a {\em
    $\Slgray$-category\/} is a category which is enriched over the
  category $\Strat$ with respect to the lax Gray tensor product
  $\laxgray$. We will adopt the notations $\ECat{\Strat\sgray}$ and
  $\ECat{\Slgray}$ to denote the corresponding (huge) categories of
  (large) enriched categories and enriched functors between them.
\end{defn}

\begin{obs}
  Of course, the lax Gray tensor $\laxgray$ is not symmetric and so
  the reader may be a little concerned that the basic theory of
  enriched categories, as expounded in Kelly~\cite{Kelly:1982:ECT}
  say, might fail to carry over in whole or part. However, it is well
  known that almost all of the fundamental definitions of enriched
  category theory may be recast in this context and that basic results
  regarding limits, colimits, Yoneda's lemma and the like have direct
  analogues in this setting. Indeed, in our case the presence of a
  dual operation $\arrow (-)\dual:\Strat->\Strat.$ which is well
  behaved with respect to $\laxgray$ (cf. observations~\ref{alt.dual}
  and~\ref{laxgray.obs} of~\cite{Verity:2005:WeakComp}) allows almost
  all of the symmetrical theory to be generalised with little
  alteration.  Although we will need no more than the (obvious) basic
  definitions here, the interested reader might like to consult
  Street's recently reprinted paper~\cite{Street:2005:EnrCatCohom} for
  an excellent introduction to the basic theory of enriched categories
  in an even more general setting, that of enrichment over a {\em
    bicategory}.
\end{obs}

\begin{defn}[complicial Gray-category]
  A {\em complicial Gray-category\/} (or usually just a Gray-category)
  is a $\Strat\sgray$-category in which each homset is a weak
  complicial set. A {\em Gray-functor\/} is defined to be a
  $\Strat\sgray$-enriched functor between Gray-categories. Let $\Gray$
  denote the (huge) category of all (possibly large) Gray-categories
  and Gray-functors between them.
\end{defn}

\begin{eg}[the Gray-category of weak complicial sets]
  \label{wcs.gray.obs.a}
  We may canonically enrich $\Strat$ with respect to its Gray tensor
  $\gray$, to obtain an enriched category denoted $\Strat\sgray$, by
  taking $\hom(X,Y)$ to be the stratified homset between the stratified
  sets $X$ and $Y$ (cf.\ Kelly~\cite{Kelly:1982:ECT} for the
  details). Now theorem~\ref{comp.hom.stab.thm}
  of~\cite{Verity:2005:WeakComp} tells us that its enriched full
  subcategory $\Wcs$ of weak complicial sets is actually a
  Gray-category in the sense of the last definition.
\end{eg}

\begin{eg}[complicial Gray-categories generalise classical ones]
  A {\em classical Gray-category\/} is simply a category enriched over
  the category $\InfCat$ of (strict) \infc-categories with respect to
  its Gray tensor product $\gray$. For example, if we restrict
  attention to to those classical Gray-categories whose homsets are
  2-categories then we obtain the kind of 3-dimensional Gray-category
  discussed by Gordon, Power and Street in~\cite{Gordon:1995:Tricats}.

  For most purposes, it is easier to present such classical
  Gray-categories in terms of the lax Gray tensor $\laxgray$ on
  $\InfCat$. Recall, from say Crans~\cite{Crans:1995:PhD} or
  Steiner~\cite{Steiner:1991:Tensor}, that if $\icat{B}$ and
  $\icat{C}$ are \infc-categories then their lax Gray tensor
  $\icat{B}\laxgray\icat{C}$ is an \infc-category generated by cells
  $\satom{\beta,\gamma}$ with $\beta \in\icat{B}$ and
  $\gamma\in\icat{C}$ subject to certain relations which ensure that
  this cell ``looks like'' a geometric product of globs.  We may
  describe the Gray tensor $\icat{B}\gray\icat{C}$ as the ``quotient''
  of the lax Gray tensor obtained by making each
  $\satom{\beta,\gamma}$ into an \infc-categorical $(r+s)$-equivalence
  whenever $\beta$ is an $r$-cell, $\gamma$ is an $s$-cell and
  $r,s>0$. It follows, therefore, that a classical Gray-category
  $\mathcal{G}$ may be presented as a category enriched over the
  category $\InfCat$ with respect to its lax Gray tensor product
  $\laxgray$ in which each composition operation $\arrow\circ:
  \mathcal{G}(b,c)\laxgray\mathcal{G}(a,b)->\mathcal{G}(a,c).$ carries
  the generators identified in the last sentence to
  \infc-categorical $(r+s)$-equivalences.

  Now theorem~255 of~\cite{Verity:2006:Complicial} tells us that
  $\arrow \inladj:\Strat->\InfCat.$, the left adjoint to Street's
  nerve functor $\arrow\inerve:\InfCat->\Strat.$, is strongly monoidal
  with respect to the {\em lax\/} Gray tensor products on those
  categories, so it follows that $\inerve$ is itself monoidal. The
  last paragraph tells us that any classical Gray-category
  $\mathcal{G}$ may be considered to be enriched with respect to
  the lax Gray tensor product, so it follows that if we apply
  $\inerve$ to the homsets of $\mathcal{G}$ we obtain a category
  enriched in (strict) complicial sets with respect to the lax Gray
  tensor $\laxgray$ on $\Strat$.

  Furthermore, it is easy to extend the calculations of loc.\ cit.\ to
  show that the Gray tensor product $\inladj(X)\gray \inladj(Y)$ may
  be constructed from the lax Gray tensor $\inladj(X)\laxgray
  \inladj(Y)\cong \inladj(X\laxgray Y)$ by making each of its
  generators of the form $\statom{(x\cdot\partproj^{r,s}_1,
    y\cdot\partproj^{r,s}_2)}$ (cf.\ observation~246 of loc.\ cit.)
  with $r,s>0$ into an \infc-categorical
  $(r+s)$-equivalence. Dualising this result using the adjunction
  $\inladj\dashv\inerve$ and the fact that $\inerve$ is full and
  faithful (cf.\ theorem 266 of loc.\ cit.), we find that the
  composition operation of $\mathcal{G}$ satisfies the condition
  described two paragraphs ago if and only if the
  compositions $\arrow\circ:\inerve(\mathcal{G}(b,c))
  \laxgray\inerve(\mathcal{G}(a,b)) -> \inerve(\mathcal{G}(a,c)).$
  carry each simplex of the form $(g\cdot\partproj^{r,s}_1,
  f\cdot\partproj^{r,s}_2)$ with $r,s>0$ to a $(r+s)$-simplex in
  $\inerve(\mathcal{G}(a,c))$ which is thin in the {\em equivalence
    stratified\/} version of Street's nerve
  $\inerve^e(\mathcal{G}(a,c))$ discussed in
  example~\ref{estrat.nerve} of~\cite{Verity:2005:WeakComp}.

  Notice, however, that the simplices $(g\cdot\partproj^{r,s}_1,
  f\cdot\partproj^{r,s}_2)$ identified in the last paragraph are
  simply those that are made thin in constructing the pre-tensor
  $\inerve(\mathcal{G}(b,c))\spretens\inerve(\mathcal{G}(a,b))$ in
  observation~\ref{pretens.tens.comp}
  of loc.\ cit. Consequently, if $\mathcal{G}$ is a
  classical Gray-category the associated composition operations of the
  last paragraph actually provide stratified maps $\arrow\circ:
  \inerve(\mathcal{G}(b,c)) \spretens\inerve(\mathcal{G}(a,b)) ->
  \inerve^e(\mathcal{G}(a,c)).$. Now, on consulting the definitions of
  the pre-tensors $\pretens$ and $\spretens$, we find that the only
  $n$-simplices which are thin in $\inerve^e(\mathcal{G}(b,c))
  \spretens\inerve^e(\mathcal{G}(a,b))$ but are not already thin in
  $\inerve(\mathcal{G}(b,c)) \spretens\inerve(\mathcal{G}(a,b))$ are
  of the form $(g\cdot\eta^n, f)$ with $f$ a thin $n$-simplex in
  $\inerve^e(\mathcal{G}(a,b)$ or $(g, f\cdot\eta^n)$ with $g$ a thin
  $n$-simplex in $\inerve^e(\mathcal{G}(b,c))$. In the first of these
  cases, we know that $g$ simply corresponds to a $0$-cell of
  $\mathcal{G}(b,c)$ and it is easily shown that the composite
  $(g\cdot\eta^n)\circ f$ in $\inerve^e(\mathcal{G}(a,c))$ is simply
  the thin $n$-simplex obtained by applying the stratified map
  $\arrow\inerve^e(g\circ-):\inerve^e(\mathcal{G}(a,b))->
  \inerve^e(\mathcal{G}(a,c)).$ to the thin $n$-simplex $f$ in
  $\inerve^e(\mathcal{G}(a,b))$. Similarly, in the second case, the
  composite $g\circ (f\cdot\eta^n)$ is simply the thin $n$-simplex
  $\inerve^e(-\circ f)(g)$, and it follows therefore that our
  composition operation actually provides us with a stratified map
  $\arrow\circ: \inerve^e(\mathcal{G}(b,c)) \spretens
  \inerve^e(\mathcal{G}(a,b)) -> \inerve^e(\mathcal{G}(a,c)).$. 

  Finally we know, from example~\ref{estrat.nerve} of loc.\ cit., that
  $\inerve^e(\mathcal{G}(a,c))$ is a weak complicial set and, from
  observation~\ref{pretens.tens.comp} of the same paper, that the
  entire inclusion $\overinc\subseteq_e:\inerve^e(\mathcal{G}(b,c))
  \spretens \inerve^e(\mathcal{G}(a,b))->\inerve^e(\mathcal{G}(b,c))
  \gray \inerve^e(\mathcal{G}(a,b)).$ is an inner anodyne
  extension. It follows that we have lifts $\arrow\circ:
  \inerve^e(\mathcal{G}(b,c)) \gray \inerve^e(\mathcal{G}(a,b)) ->
  \inerve^e(\mathcal{G}(a,c)).$ providing composition operations which
  make the weak complicial sets $\inerve^e(\mathcal{G}(a,b))$ into the
  homsets of a complicial Gray-category. Finally it is clear that this
  construction faithfully represents classical Gray-categories as
  complicial ones.
\end{eg}

\begin{obs}[size]
  In what follows, we shall work somewhat informally with categorical
  structures that naturally live in different {\em universes\/} of
  sets. To formalise things the reader might like to think in terms of
  three fixed and nested Grothendieck universes of sets, whose
  denizens are called {\em small}, {\em large\/} and {\em huge\/} sets
  respectively. For example, in~\cite{Verity:2005:WeakComp} our
  stratified sets were generally drawn from the world of small sets
  and so the categories $\Strat$ and $\Wcs$ were large and locally
  small (small homsets). Our new category $\Gray$ is a huge structure.

  This hierarchy all works very well, until we wish to start talking
  about the nerves of Gray-categories. These latter structures are
  large so their nerves should also be large structures and thus are
  not strictly speaking members of the category $\Strat$ of small
  stratified sets. One solution to this problem might be to adopt a
  separate notation $\category{STRAT}$ for the category of large
  stratified sets, which can then act as a carrier for the nerves of
  large Gray-categories. This, however, simply serves to complicate
  our notation and confuse matters. 

  We instead take the approach of overloading notations such as
  $\Strat$ to denote a structure inhabiting whichever universe is
  necessary for the current argument. This then allows us to write
  $\arrow \nerv:\Gray->\Strat.$ for the nerve construction we will discuss
  and to infer that in this case $\Strat$ denotes the category of
  large stratified sets.
\end{obs}

\section{The Classical Theory of Homotopy Coherent Nerves}

Now we recall (a suitable version of) the homotopy coherent nerve
construction for categories enriched in simplicial sets (called {\em
  simplicially enriched categories}). Much of the work presented in
this section is due to Cordier~\cite{Cordier:1982:HtyCoh} and Cordier
and Porter~\cite{Cordier:1986:HtyCoh}.

\begin{defn}[the homotopy coherent \infc-path]\label{hcpath.defn}
  We define a locally ordered category $\pspath$ whose objects are
  the integers and whose homsets are {\em ordered cubes}, that is to
  say that they are all suitable powers 
  \begin{equation*}
    [1]^n\defeq\underbrace{[1]\times[1]\times\cdots\times[1]}_{n\text{ factors}}
  \end{equation*}
  of the two point ordinal $[1]$. We recall briefly that, as an
  iterated product of ordered sets, such powers are ordered
  ``pointwise'', that is to say $(a_n,a_{n-1},...,a_1)\leq
  (b_n,b_{n-1},...,b_1)$ iff $a_i\leq b_i$ for all
  $i=1,2,...,n$. We choose to index the ordinates of our
  tuples from right to left, a convention we adopt in order to
  simplify our notation later on.

  For integers $r<s$ our homset $\pspath(r,s)$ will be defined to be
  isomorphic to $[1]^{(s-r-1)}$, however for reasons of combinatorial
  convenience we actually elect to specify it as a subset of
  $[1]^{(s-t)}$:
  \begin{equation*}
    \pspath(r,s) =
    \begin{cases}
      \emptyset & \text{if $s<r$,} \\
      \{ () \} = [1]^0 & \text{if $s=r$,} \\ 
      \{(a_s,a_{s-1},...,a_{r+1})\in[1]^{(s-r)}\mid a_s=0\} &
      \text{if $s>r$.}
    \end{cases}
  \end{equation*}
  A consequence of this convention is that while $\pspath(r,r)$ and
  $\pspath(r,r+1)$ are both isomorphic to the single point set, their
  respective elements are clearly distinguished in intention. Notice
  also that we adopt the convention of indexing the elements of an
  arrow in $\pspath(r,s)$ using the integers $r+1,r+2,...,s$, which
  further simplifies and illuminates many of the arguments we make
  later on.

  The composition of $\pspath$ is now a matter of mere concatenation,
  in other words if $\vec{a} = (a_s,a_{s-1},...,a_{r+1})$ is an
  element of $\pspath(r,s)$ and $\vec{b}=(b_t,b_{t-1},...,b_{s+1})$ is
  an element of $\pspath(s,t)$ then their composite is defined by
  \begin{equation*}
    (b_t,b_{t-1},...,b_{s+1})\circ(a_s,a_{s-1},...,a_{r+1}) \defeq
    (b_t,...,b_{s+1},a_s,...,a_{r+1})
  \end{equation*}
  which is clearly an element of $\pspath(r,t)$. Of course this
  operation is associative and it preserves the orderings of the
  homsets since they are defined pointwise. It also immediately
  illuminates why we defined $\pspath(r,r)$ to contain only the unique
  $0$-tuple $()$, since this is the identity for concatenation of
  tuples.  It follows therefore that $\pspath$ is a well-defined
  locally ordered category.

  Our real interest, however, is in defining an associated
  simplicially enriched category $\hcpath$ by applying the classical
  categorical nerve functor $\arrow \nerv:\Cat->\Simp.$ (cf.\
  observation~\ref{cat.nerves} of~\cite{Verity:2005:WeakComp}) to each
  of the homsets of $\pspath$. Abstractly we find that this process
  gives us a category enriched over $\Simp$ simply because we know
  that $\nerv$ preserves products (it is a right adjoint).  More
  concretely, however, the $l$-simplices of $\hcpath(r,s)$ are ordered
  sequences of tuples $\vec{a}_0\leq \vec{a}_1\leq\cdots\leq\vec{a}_l$
  of length $(l+1)$ in $\pspath(r,s)$ and these are composed
  pointwise, in other words if $\vec{b}_0\leq
  \vec{b}_1\leq\cdots\leq\vec{b}_l$ is a $l$-simplex of $\hcpath(s,t)$
  then their composite is $\vec{b}_0\circ\vec{a}_0\leq
  \vec{b}_1\circ\vec{a}_1\leq\cdots \leq\vec{b}_l\circ\vec{a}_l$.

  We might call $\hcpath$ the {\em generic homotopy coherent
    \infc-path}, since simplicial functors out of it and into any
  other simplicially enriched category $\mathcal{E}$ correspond to
  weakened functors from the ordered set $\mathbb{N}$ into
  $\mathcal{E}$ which {\em preserve composition up to coherent
    homotopy}. It was this point of view which originally motivated
  Cordier to study structures of this kind~\cite{Cordier:1982:HtyCoh}.

  Before moving on we will also adopt the notation $\hcpath[n]$
  (resp.\ $\pspath[n]$) for the (enriched) full subcategory of
  $\hcpath$ (resp.\ $\pspath$) whose objects are the elements of the
  ordinal $[n]$, which we might think of as being the {\em generic
    homotopy coherent $n$-path}.
\end{defn}

\begin{obs}[$\pspath$ as a free structure]\label{pspath.free}
  We may also describe $\pspath$ as a freely generated locally ordered
  category. To do so, observe that if we have an arrow
  $\vec{a}=(a_s,...,a_{r+1})$ in $\pspath(r,s)$ then we may find some
  index $r<k<s$ for which $a_k=0$ if and only if we may decompose
  $\vec{a}$ as a composite of two non-identity arrows
  $(a_k,...,a_{r+1})$ in $\pspath(r,k)$ and $(a_s,...,a_{k+1})$ in
  $\pspath(k,s)$. Consequently there is precisely one arrow in
  $\pspath(r,s)$ which is not decomposable in this way, that being the
  one consisting of a sequence of $(s-r-1)$ copies of $1$ followed by a single
  $0$ and for which we adopt the denotation $\indec<r,s>$. Furthermore any
  $\vec{a}\in\pspath(r,s)$ may be ``split at its $0$s'' to express it
  uniquely as a composite of indecomposable arrows. More precisely if
  $r_1<r_2<\cdots<r_k=s$ enumerate the indices at which the $\vec{a}$
  has a ordinate with value $0$ then we may uniquely express it as the
  composite $\indec<r_{k-1},r_k>\circ\cdots\circ\indec<r_1,r_2>\circ
  \indec<r,r_1>$.

  This certainly demonstrates that $\pspath$ is freely generated as a
  mere category by its indecomposable arrows, indeed we've shown that
  it is the free category associated with the (reflexive) graph
  underlying the ordered set $\mathbb{N}$ (qua category). We can go
  further, however, and describe its local orders in a similar
  fashion, by observing that if $\vec{a}\leq\vec{b}$ in $\pspath(r,s)$
  then $\vec{a}$ has a $0$ wherever $\vec{b}$ has one and thus that
  the splitting process described above (at least) splits the former
  wherever it splits the latter. Consequently any such inequality may
  be obtained by using $\circ$ to compose inequalities of the form
  $\indec<r_{k-1},r_k>\circ\cdots\circ\indec<r_1,r_2>\circ
  \indec<r_0,r_1>\leq \indec<r_0,r_k>$ which in turn may be expressed
  canonically in terms of the most primitive inequalities of the form
  $\indec<s,t>\circ\indec<r,s>\leq\indec<r,t>$ (for $r<s<t$).

  Summarising all of this we see that, as a locally ordered category,
  $\pspath$ is freely generated by its indecomposable arrows
  $\indec<r,s>$ subject to the primitive inequalities
  $\indec<s,t>\circ\indec<r,s>\leq\indec<r,t>$. In other words, if
  $\lcat{E}$ is any other locally ordered category then to (uniquely)
  define a locally ordered functor $\arrow f:\pspath->\lcat{E}.$ it is
  enough to specify its action $f(r)$ on objects and $f\indec<r,s>$ on
  indecomposable arrows and to check that $f(r)=\dom(f\indec<r,s>)$ and
  $f(s)=\cod(f\indec<r,s>)$ in $\lcat{E}$ (for all $r<s$) and that
  $f\indec<s,t>\circ f\indec<r,s>\leq f\indec<r,t>$ in the ordered set
  $\lcat{E}(f(r),f(t))$ (for all $r<s<t$).
\end{obs}

\begin{obs}\label{hcpath.functor}
  We may now construct a functor from $\Delta$ to the category of
  locally ordered categories, which carries $[n]$ to $\pspath[n]$. To
  define its action on simplicial operators $\arrow \alpha:[n]->[m].$
  we see from the last observation that it is enough to specify how
  the order enriched functor $\arrow
  \pspath(\alpha):\pspath[n]->\pspath[m].$ should act on
  indecomposable arrows $\indec<r,s>$ in $\pspath[m]$. This, however,
  is easily done by letting 
  \begin{equation}\label{pspath.op.act}
    \pspath(\alpha)\indec<r,s>\defeq\indec<\alpha(r), \alpha(s)>
  \end{equation}
  under the convention that the notation $\indec<t,t>$, which becomes
  an issue if $\alpha(r)=\alpha(s)$, is taken to mean the identity
  $0$-tuple $()$ in $\pspath(t,t)$. The verifications of the
  conditions given in the final paragraph of the last observation are
  now a trivial matter, thus giving the required functor. Furthermore,
  if $\arrow\beta:[m]->[p].$ is a second simplicial operator than we
  may show that the functors $\pspath(\beta)\circ\pspath(\alpha)$ and
  $\pspath(\beta\circ\alpha)$ are identical by checking that they
  coincide on indecomposables, simply using the definition
  in~(\ref{pspath.op.act}), and then appealing to the uniqueness
  clause of the last observation.

  Again, our intent here is actually to use the functor of the last
  paragraph to construct a functor $\hcpath$ from $\Delta$ to the
  category $\ECat\Simp$ of simplicially enriched categories and
  simplicial functors. The we may now do simply by applying the
  categorical nerve construction $\arrow \nerv:\Cat->\Simp.$ to the
  homsets and maps of homsets that arise when we apply the functor
  $\pspath$. Explicitly, this functor maps $[n]$ to the simplicially
  enriched category $\hcpath[n]$ of observation~\ref{hcpath.defn} and
  it maps a simplicial operator $\arrow\alpha:[n]->[m].$ to a
  simplicial functor $\arrow\hcpath(\alpha):\hcpath[n]->\hcpath[m].$
  which acts like $\alpha$ on objects and applies $\pspath(\alpha)$
  pointwise to each $k$-simplex $\vec{a}_0\leq
  \vec{a}_1\leq\cdots\leq\vec{a}_k$ in $\hcpath(r,s)=\hcpath[n](r,s)$
  to give $\pspath(\alpha)(\vec{a}_0)\leq
  \pspath(\alpha)(\vec{a}_1)\leq\cdots \leq\pspath(\alpha)(\vec{a}_k)$
  in $\hcpath(\alpha(r),\alpha(s))=\hcpath[m]( \alpha(r),\alpha(s))$.
\end{obs}

\begin{recall}[homotopy coherent nerve]\label{hcnerv.simp.def}
  We may apply Kan's construction \cite{Kan:1958:Adjoint} to the
  functor $\hcpath$ of the last observation and form an adjoint pair of
  functors:
  \begin{equation*}
    \adjdisplay \ECat\Simp -> \Simp, \ladj_{hc}-| \nerv_{hc}.
  \end{equation*}
  The functor $\nerv_{hc}$ is customarily called the {\em homotopy
    coherent nerve construction}. It carries a simplicially enriched
  category $\lcat{E}$ to a simplicial set whose $n$-simplices are
  simplicial functors $\arrow f:\hcpath[n]->\lcat{E}.$ and whose
  simplicial action is given by pre-composition $f\cdot\alpha\defeq
  f\circ \hcpath(\alpha)$. Furthermore, its left adjoint may be
  expressed as the colimit $\ladj_{hc}(X)\cong\colim(X,\hcpath)$ of
  the diagram $\hcpath$ in $\ECat\Simp$ {\em weighted by\/} the
  presheaf $\arrow X:\Delta\op->\Set.$
  (cf.\ Kelly~\cite{Kelly:1982:ECT}).
\end{recall}

\begin{obs}[horns in homotopy coherent nerves]\label{hchorn.obs}
  In~\cite{Cordier:1986:HtyCoh} Cordier and Porter demonstrate that
  the homotopy coherent nerve of a simplicially enriched category
  whose homsets are all Kan complexes is actually a quasi-category
  (weak Kan complex). We will prove a generalisation of this result
  for Gray-categories below.  For now however it will be useful to
  provide an explicit description of the inner horns $\arrow
  h:\Lambda^k[n]->\nerv_{hc}(\lcat{E}).$ of a homotopy coherent
  nerve. Applying the adjunction $\ladj_{hc}\dashv \nerv_{hc}$ we see that
  these correspond to simplicial functors
  $\arrow\hat{h}:\ladj_{hc}(\Lambda^k[n])->\lcat{E}.$ so this task
  immediately reduces to that of providing an explicit description of
  the simplicially enriched category $\ladj_{hc}(\Lambda^k[n])$. Indeed,
  we will actually show that this bears a natural description as a
  subcategory of $\ladj_{hc}(\Delta[n])\cong\hcpath[n]$.
\end{obs}


\begin{obs}[$\hcpath$ as a free structure]\label{Fhc.freeness}
  Forgetting, for the moment, the simplicial set structures on the
  homsets of $\hcpath$ we obtain a category which may again be
  described as the free category on some reflexive graph.  In
  particular, following the pattern laid down in
  observation~\ref{pspath.free}, we consider an {\em $l$-arrow\/}
  $\vec{a}_0\leq\vec{a}_1\leq\cdots\leq\vec{a}_l$ of $\hcpath$ and
  observe that if its topmost tuple $\vec{a}_l$ has a $0$ at some
  index then each $\vec{a}_i$ must also have a $0$ at that index,
  since they are each less than or equal to the former in the
  pointwise ordering. It follows that each tuple in
  $\vec{a}_0\leq\vec{a}_1\leq\cdots\leq\vec{a}_l$ ``may be split at
  the $0$s of $\vec{a}_l$'' to express this $l$-arrow as a unique
  (pointwise) composite of indecomposable arrows, these being those
  arrows of $\hcpath$ whose topmost tuple is an indecomposable arrow
  $\indec<r,s>$ of $\pspath$.  Furthermore, if
  $\arrow\alpha:[n]->[m].$ is a simplicial operator then by definition
  the simplicial functor
  $\arrow\hcpath(\alpha):\hcpath[n]->\hcpath[m].$ carries an
  indecomposable $l$-arrow with topmost tuple $\indec<r,s>$ to an
  $l$-arrow with topmost tuple $\indec<\alpha(r),\alpha(s)>$, which is
  thus either indecomposable itself or an identity (when
  $\alpha(r)=\alpha(s)$). Consequently, if we regard $\hcpath$ as a
  functor from $\Delta$ to $\Cat$ then we may factor it (up to
  isomorphism) through the free category functor $\arrow
  \ladj:\Grph->\Cat.$, from reflexive graphs to categories, by restricting
  the image $\hcpath[n]$ of each object $[n]\in\Delta$ to its
  subgraph of indecomposables to obtain a functor
  $\arrow\hcpathg:\Delta->\Grph.$.

  The utility of this observation lies in the fact that it enables us
  to evaluate the colimits used to define $\ladj_{hc}$ in
  observation~\ref{hcnerv.simp.def}. First we may happily calculate
  without worrying about the simplicial structure of homsets, since
  the forgetful functor $\overarr:\ECat\Simp->\Cat.$ preserves and
  reflects (indeed creates) the colimits of $\ECat\Simp$.  Then we may
  construct the colimit $\ladj_{hc}\cong\colim(X,\hcpath)$ in $\Cat$ by
  evaluating the corresponding colimit $\colim(X,\hcpathg)$ in
  $\Grph$, where such things are constructed as in $\Set$, and then
  applying the free category functor. 
\end{obs}



\begin{obs}\label{ident.hc.horn}
  For instance, the last observation allows us to show that $\ladj_{hc}$
  carries inclusions of simplicial sets to inclusions of simplicially
  enriched categories. To do so first observe that the functor
  $\arrow\pspath(\face_i):\pspath[n-1]-> \pspath[n].$ of locally
  ordered categories is injective and that its image contains the
  arrow $\vec{a}\in\pspath(r,s)$ if and only if $s<i$ or $i<r$ or
  ($r<i<s$ and $a_i=1$).  This latter fact may be ascertained simply
  by checking that it holds for indecomposables $\indec<r,s>$ and
  extending to all arrows using the decomposition at $0$s result. It
  follows immediately that
  $\arrow\hcpathg(\face_i):\hcpathg[n-1]->\hcpathg[n].$ in $\Grph$ is
  also injective and that the indecomposable $l$-arrow
  $\vec{a}_0\leq\vec{a}_1\leq\cdots \leq\vec{a}_l$ is in it image iff
  its bottommost tuple $\vec{a}_0$, and thus each one of its members,
  satisfies the condition of the last sentence. Applying this
  characterisation it is also clear that for each $i<j$ we have a
  pullback
  \begin{equation}\label{pb.diamond}
    \xymatrix@R=0.75em@C=3em{
      &
      {\hcpathg[n-2]}\ar@{u(->}[rd]^-{\hcpathg(\face_{j-1})}
      \ar@{u(->}[ld]_-{\hcpathg(\face_i)}
      & \\
      {\hcpathg[n-1]}\ar@{u(->}[dr]_-{\hcpathg(\face_j)} &&
      {\hcpathg[n-1]}\ar@{u(->}[ld]^-{\hcpathg(\face_i)} \\
      & {\hcpathg[n]} & }
  \end{equation}
  in $\Grph$. So consider the ``wide pushout'' diagram in simplicial
  sets consisting of $(n+1)$ copies of $\Delta[n-1]$ indexed by
  elements of $[n]$ and $(n+1)n/2$ copies of $\Delta[n-2]$ indexed by
  pairs $i,j\in[n]$ with $i<j$ each of which comes equipped with a
  pair of simplicial maps into two of the copies of $\Delta[n-1]$ as
  depicted below.
  \begin{equation}\label{pb.diamond.2}
    \xymatrix@R=0.75em@C=3em{
      & {\Delta[n-2]_{i,j}}\ar@{u(->}[rd]^-{\Delta(\face_{j-1})}
      \ar@{u(->}[ld]_-{\Delta(\face_i)} & \\
      {\Delta[n-1]_j}\ar@{u(..>}[dr]_-{\Delta(\face_j)} && 
      {\Delta[n-1]_i}\ar@{u(..>}[ld]^-{\Delta(\face_i)} \\
      & {\Delta[n]} & }
  \end{equation}
  The dotted inclusions in this picture provide a cone under our
  diagram and each of the diamonds above is a pullback in
  $\Simp$. Therefore, since the colimits of $\Simp$ are calculated as
  in $\Set$, it is easily seen that the colimit of our wide pushout
  diagram is isomorphic to the simplicial set obtained by taking the
  union of the images of the dotted inclusions in $\Delta[n]$, which
  is simply the boundary $\boundary\Delta[n]$.

  Now recall that colimits weighted by representables are easily
  calculated by ``evaluation'', in other words we have a natural
  isomorphism $\colim(\Delta[n],\hcpathg)\cong\hcpathg[n]$. It follows
  that if we apply the functor $\colim(-,\hcpathg)$ to the diamond in
  display~(\ref{pb.diamond.2}) then we obtain the one in
  display~(\ref{pb.diamond}). Furthermore the weighted colimit
  construction preserves colimits in each variable, so in particular
  the functor $\colim(-,\hcpathg)$ preserves the colimit of our wide
  pushout $\boundary\Delta[n]$. However, arguing just as we did
  before, using the fact that colimits in $\Grph$ are calculated as in
  $\Set$, we see that the pullbacks of display~(\ref{pb.diamond})
  allow us to show that the colimit of this wide pushout in $\Grph$
  is also isomorphic to the subgraph of $\hcpathg[n]$ obtained by
  taking the union of the images of the inclusions
  $\inc\hcpathg(\face_i):\hcpathg[n-1]-> \hcpathg[n].$. 
  In other words, we have demonstrated that the induced map
  \begin{equation}\label{colim.bdary.inc}
    \arrow\colim(\subseteq_s,\hcpathg):\colim(\boundary\Delta[n],\hcpathg)->
    \colim(\Delta[n],\hcpathg).
  \end{equation}
  restricts to an isomorphism between its codomain and the subgraph
  of $\hcpathg[n]\cong\colim(\Delta[n],\hcpathg)$ identified in the
  last sentence.  

  However, we also know that every inclusion $\inc i:X->Y.$ of
  simplicial sets may be constructed as a transfinite composite of
  pushouts of the boundary inclusions
  $\inc\subseteq_s:\boundary\Delta[n]->\Delta[n].$. So applying the
  colimit preservation property of $\colim(-,\hcpathg)$ again we see
  that the induced map $\inc\colim(i,\hcpathg):
  \colim(X,\hcpathg)->\colim(Y,\hcpathg).$ may be constructed as a
  transfinite composite of pushouts of the inclusions of
  display~(\ref{colim.bdary.inc}), and is thus itself an inclusion
  since these operations clearly preserve inclusions in $\Set$ and
  therefore also do so in $\Grph$. Finally, applying the free category
  functor, which carries inclusions of graphs to inclusions of
  categories, and applying observation~\ref{Fhc.freeness} we
  find that the simplicial functor $\arrow
  \ladj_{hc}(i):\ladj_{hc}(X)->\ladj_{hc}(Y).$ is also an inclusion of
  simplicially enriched categories. In particular, the simplicial
  functor
  \begin{equation}\label{horn.inc.scat}
    \arrow \ladj_{hc}(\subseteq_s):\ladj_{hc}(\Lambda^k[n])->\ladj_{hc}(\Delta[n])
    \cong\hcpath[n].
  \end{equation}
  is an inclusion, as suggested above, and its image $\hchorn^k[n]$,
  called {\em the homotopy coherent $k$-horn}, can be characterised as
  being the smallest subcategory of $\hcpath[n]$ containing the images
  of $\inc \hcpath(\face_i):\hcpath[n-1]->\hcpath[n].$ for each $i\in
  [n]\setminus\{k\}$.
\end{obs}


\section{Nerves of $\protect\Strat\protect\slgray$-Categories}

We seek to extend this classical material and apply it to
the problem of providing a well behaved nerve construction which
carries Gray-categories to weak complicial sets. To do so we first
provide the homsets of $\hcpath$ with an appropriate stratification.

\begin{obs}\label{redesc.hcpath.obs}
  It will be convenient for what follows to slightly rephrase the
  definition of $\hcpath$ given in definition~\ref{hcpath.defn}. To do
  so recall that the categorical nerve construction $\arrow
  \nerv:\Cat->\Simp.$ preserves products and carries each ordinal $[n]$ to
  the standard simplex $\Delta[n]$. Furthermore, we know that the
  homset $\pspath(r,s)$ ($r\leq s$) is defined to be a subset of the
  iterated product $[1]^{(s-r)}$, so it follows that we can think of
  $\hcpath(r,s)$, its nerve, as the simplicial subset of
  $\nerv([1]^{(s-r)})\cong \Delta[1]^{(s-r)}$ of those simplices
  $\vec\alpha=(\alpha_s,...,\alpha_{r+1})$ for which $\alpha_s$ is the
  constant $0$ operator defined in notation~\ref{simp.of.1.a}.  Under
  this representation the composition operation is again given by
  concatenation of tuples of operators. More abstractly, it will be
  useful to re-cast this by saying that composition is a restriction
  of the canonical associativity isomorphism of the monoidal category
  $(\Simp,\times,\Delta[0])$ depicted in the following square:
  \begin{equation}\label{comp.as.rest}
    \xymatrix@=2em{
      {\hcpath(s,t)\times\hcpath(r,s)}
      \ar[r]^-{\circ}\ar@{u(->}[d]_{\subseteq_s} &
      {\hcpath(r,t)}
      \ar@{u(->}[d]^{\subseteq_s} \\
      {\Delta[1]^{(t-s)}\times\Delta[1]^{(s-r)}}
      \ar[r]_-{\cong} & {\Delta[1]^{(t-r)}}}
  \end{equation}
  Of course, it is immediate that an arrow $\vec\alpha$ in
  $\hcpath(r,s)$ is decomposable iff there exists some integer $k$
  with $r<k<s$ and for which $\alpha_k$ is equal to the constant
  operator $0$ of notation~\ref{simp.of.1.a}.

  We will also find it useful to identify, in these terms, the arrows
  of $\hcpath[n]$ which are in the {\em inner\/} homotopy coherent
  horn $\hchorn^k[n]$ of observation~\ref{ident.hc.horn}. Recasting
  the characterisation derived in the first paragraph of that
  observation, we find that our arrow $\vec\alpha$ of
  $\hcpath[n](r,s)$ is in the image of a face inclusion
  $\inc\hcpath(\face_i):\hcpath[n-1]->\hcpath[n].$ iff $s<i$ or $i<r$
  or ($r<i<s$ and $\alpha_i=1$). In the case of $\hcpath(\face_0)$
  this is simply the full subcategory of $\hcpath[n]$ on the objects
  $1,2,...,n$, whereas for $\hcpath(\face_n)$ it is the full
  subcategory on objects $0,1,...,n-1$. This immediately implies that
  $\hchorn^k[n]$ completely contains all homsets of $\hcpath[n]$
  except for $\hcpath[n](0,n)$. Furthermore, any decomposable arrow of
  $\hcpath[n](0,n)$ factors into a composite of arrows in
  $\hcpath[n](0,r)$ and $\hcpath[n](r,n)$ for some $0<r<n$, but the
  former homset is in the image of $\hcpath(\face_n)$ and the latter
  is in the image of $\hcpath(\face_0)$ so it follows that their
  composite is in the subcategory generated by the union of these
  images. Combining this observation with our characterisation of
  decomposable arrows and throwing in the images of $\hcpath(\face_i)$
  for the remaining faces with $i\neq k$, we find that the arrow
  $\vec\alpha$ of $\hcpath(0,n)$ is in $\hchorn^k[n](0,n)$ iff there
  exists some integer $i$ with $0<i<n$ for which either $\alpha_i=0$
  or ($i\neq k$ and $\alpha_i=1$).

  After a few moments reflection, it becomes clear that this latter
  characterisation may itself be re-expressed in terms of the {\em
    corner product\/} $\ctimes$ associated with the product of
  simplicial sets as defined in observation~\ref{corner.tensor}
  of~\cite{Verity:2005:WeakComp}. To be precise, the simplicial map
  which ``forgets the leftmost ordinate'' provides us with an
  isomorphism between $\hcpath[n](0,n)$ and $\Delta[1]^{(n-1)}$ the
  codomain of the following iterated corner product
  \begin{equation}\label{horn.corner.prod}
    (\overinc\subseteq_r:\boundary\Delta[1]->\Delta[1].)^{\ctimes(n-1-k)}
    \ctimes(\overinc\subseteq_r:\Lambda^1[1]->\Delta[1].)\ctimes
    (\overinc\subseteq_r:\boundary\Delta[1]->\Delta[1].)^{\ctimes(k-1)}
  \end{equation}
  and this restricts to provide us with an isomorphism between its
  domain and $\hchorn^k[n](0,n)$.
\end{obs}

\begin{obs}[the nerve of homotopy coherent paths]\label{hcpath.of.graycat}
  Consider the functor $\arrow\triv_0:\Simp->\Strat.$ which stratifies
  each simplicial set using the maximal stratification. This is right
  adjoint to the stratification forgetting functor and thus preserves
  all limits, so in particular it provides us with a strict monoidal
  functor from $(\Simp,\times,\Delta[0])$ to $(\Strat,\gray, \Delta[0]
  )$. Consequently, we may apply this trivialisation locally to the
  homsets of simplicially enriched categories to obtain a functor
  $\arrow\ECat{\triv_0}:\ECat\Simp->\ECat{\Strat\sgray}.$. Composing this
  with the functor $\arrow\hcpath:\Delta->\ECat\Simp.$ of
  observation~\ref{hcpath.functor} we obtain a functor which we shall
  also denote by $\hcpath$ and from which we may derive an adjunction
  \begin{equation*}
     \adjdisplay {\ECat{\Strat\sgray}} -> \Simp, \ladj_{hc}-| \nerv_{hc}.
  \end{equation*}
  by Kan's construction. The simplicial set $H_{hc}(\lcat{E})$ is
  called the {\em nerve of homotopy coherent paths}.


  Now suppose that $\lcat{E}$ is a $\Strat\sgray$-category, then
  applying the analysis of horns in homotopy coherent nerves, which we
  commenced in observation~\ref{hchorn.obs}, it becomes clear that the
  nerve $\nerv_{hc}(\lcat{E})$ is a quasi-category iff $\lcat{E}$ has the
  RLP with respect to each inner homotopy coherent horn inclusion
  $\overinc: \hchorn^k[n]->\hcpath[n].$, where again $\hchorn^k[n]$ is
  stratified by applying $\triv_0$ to its homsets. In turn, this
  reduces to showing that each homset $\lcat{E}(e,e')$ of $\lcat{E}$
  has the RLP with respect to
  $\overinc\subseteq_r:\hchorn^k[n](0,n)->\hcpath[n](0,n).$, the only
  homset inclusion at which $\hchorn^k[n]$ and $\hcpath[n]$
  differ. Returning to observation~\ref{redesc.hcpath.obs} we see that
  this latter inclusion is isomorphic to the one obtained by applying
  $\triv_0$ to the the corner product in
  display~(\ref{horn.corner.prod}) and appealing to the monoidal
  properties of $\triv_0$ we see that this is equal to the iterated
  corner tensor
  \begin{equation*}
    (\overinc\subseteq_r:\boundary\Delta[1]->\Delta[1]_t.)^{\cgray(n-1-k)}
    \cgray(\overinc\subseteq_r:\Lambda^1[1]->\Delta^1[1].)\cgray
    (\overinc\subseteq_r:\boundary\Delta[1]->\Delta[1]_t.)^{\cgray(k-1)}
  \end{equation*}
  in $\Strat$. Finally, the inner factor
  $\overinc\subseteq_r:\Lambda^1[1]->\Delta^1[1].$ of this corner
  tensor is an elementary anodyne extension and its remaining factors
  are all inclusions so we may apply
  corollary~\ref{anod.tensor.stab.cor} of~\cite{Verity:2005:WeakComp}
  to show that this inclusion enjoys the LLP with respect to all weak
  complicial sets, and thus in particular that it does so with respect
  to the homsets of all Gray-categories. It follows therefore that the
  nerve $\nerv_{hc}(\lcat{G})$ is a quasi-category whenever $\lcat{G}$ is
  a Gray-category.
\end{obs}

\begin{obs}
  The result of the last observation is somewhat uninspiring, since
  all we have done is to reprove Cordier and Porter's analysis of
  homotopy coherent nerves~\cite{Cordier:1986:HtyCoh}. It is important
  to note that the nerve $\nerv_{hc}(\lcat{G})$ takes no account
  whatsoever of the non-thin simplices in the homsets of
  $\lcat{G}$. To rectify this deficiency we propose to stratify
  $\hcpath$ more carefully as a $\Slgray$-category, by replacing the
  use of cartesian products of simplicial sets in
  observation~\ref{redesc.hcpath.obs} by lax Gray tensor products of
  the corresponding stratified sets and insisting that $\hcpath(r,s)$
  is given a stratification which makes it into a regular subset of the
  $(s-r)$-fold tensor power of $\Delta[1]$.  By doing so in the next
  few pages we will construct a stratified nerve which {\em
    faithfully\/} represents Gray-categories as weak complicial sets.
\end{obs}

\begin{defn}
  Let $\hclgray$ denote the $\Slgray$-category whose
  underlying simplicially enriched category, obtained by forgetting
  stratifications on homsets, is $\hcpath$ and for which each homset
  $\hclgray(r,s)$ ($r\leq s$) is a regular subset of the
  iterated tensor:
  \begin{equation*}
    \Delta[1]^{\laxgray (s-r)}\defeq\underbrace{\Delta[1]\laxgray[1]
        \laxgray\cdots\laxgray\Delta[1]}_{(s-r)\text{ factors}}
  \end{equation*}
\end{defn}

\begin{obs}
  This definition depends upon observation~\ref{laxgray.obs}
  of~\cite{Verity:2005:WeakComp}, wherein we find that the monoidal
  structure associated with the tensor $\laxgray$ is completely
  determined by the fact that the stratification forgetting functor is
  a strict monoidal functor from $(\Strat,\laxgray,\Delta[0])$ to
  $(\Simp,\times, \Delta[0])$. First of all, this ensures that if we
  forget the stratifications of the homsets of any $\Slgray$-category
  we do actually obtain a genuine simplicially enriched category. It
  also ensures that $\hclgray$ is completely determined by the given
  conditions, since the stratified maps that make up its structure are
  themselves determined by their underlying simplicial maps which
  we've specified must be the corresponding components of the
  structure of $\hcpath$. Indeed, in order to demonstrate the
  existence of $\hclgray$ all we need do is show that the compositions
  of $\hcpath$ respect the stratifications we've specified for our
  homsets. To that end, consider the stratified isomorphism at the
  bottom of the following square
  \begin{equation}
    \xymatrix@=2em{
      {\hclgray(s,t)\laxgray\hclgray(r,s)}
      \ar[r]^-{\circ}\ar@{u(->}[d]_{\subseteq_r} &
      {\hclgray(r,t)}
      \ar@{u(->}[d]^{\subseteq_r} \\
      {\Delta[1]^{\laxgray(t-s)}\laxgray\Delta[1]^{\laxgray(s-r)}}
      \ar[r]_-{\cong} & {\Delta[1]^{\laxgray(t-r)}}}
  \end{equation}
  obtained from the monoidal structure associated with
  $\laxgray$. Since the stratification functor is strictly monoidal we
  know that the simplicial map underlying this isomorphism is simply
  the corresponding isomorphism for $\times$ as depicted in
  display~(\ref{comp.as.rest}). So we may define the composition of
  $\hclgray$ to be the restriction to {\em regular\/} subsets
  depicted and immediately infer that its underlying simplicial map is
  the composition of $\hcpath$ as postulated.

  Notice that the lax Gray tensor $\laxgray$ preserves regularity, in
  the sense that if $f$ and $g$ are regular stratified maps then so is
  $f\laxgray g$ (see observation~132
  of~\cite{Verity:2006:Complicial}). It follows that the map
  $\inc\Delta(\vertex^1_0)\laxgray \Delta[1]^{\laxgray(s-r-1)} :
  \Delta[1]^{\laxgray(s-r-1)} ->\Delta[1]^{\laxgray(s-r)}.$ is a
  regular inclusion and that it restricts to a canonical isomorphism
  between $\hclgray(r,s)$ and $\Delta[1]^{\laxgray(s-r-1)}$.
\end{obs}

\begin{obs}\label{hcpath.lift.obs}
  To show that the functor $\arrow \hcpath:\Delta->\ECat\Simp.$ also
  lifts to a functor $\arrow\hclgray:\Delta->\ECat{
    \Slgray}.$ we again apply the methodology of the last
  observation. In other words, we first show that we may describe some
  aspect of the structure of $\hcpath$ abstractly in terms of the
  monoidal structure of $(\Simp,\times,\Delta[0])$, then make the
  analogous construction for $\hclgray$ using the monoidal
  structure of $(\Strat,\laxgray,\Delta[0])$ and finally appeal to the
  to the strong monoidality of the forgetful functor to show that the
  former simplicial structure underlies the newly constructed
  stratified structure.

  To simplify matters, observe that it is enough to analyse the
  simplicial functors obtained by applying $\hcpath$ to each of the
  elementary face and degeneracy operators, since they generate
  $\Delta$. So consider these in turn:
  \begin{itemize}[fullwidth, leftmargin=1em, itemsep=1ex]
  \item {\boldmath $\arrow
      \hcpath(\face_k):\hcpath[n-1]->\hcpath[n].$} This acts on
    objects $r\in[n-1]$ simply by applying $\face_k$ and we may
    ascertain its action on homsets by first considering how the
    corresponding functor $\pspath(\face_k)$ acts on an indecomposable
    $\indec<r,s>$, as we did in observation~\ref{ident.hc.horn}. On
    doing so we find that
    \begin{equation*}
      \pspath(\face_k)\indec<r,s> =
      \begin{cases}
        \indec<r,s> & \text{if $s<k$,} \\
        \indec<r+1,s+1> & \text{if $k\leq r$,} \\
        \indec<r,s+1> & \text{if $r<k\leq s$.}
      \end{cases}
    \end{equation*}
    where we might think of the second case as shifting our
    indecomposable left to make space for a new symbol at index $k$
    below it and the third case as the actual insertion of a new $1$
    symbol at index $k$ in our indecomposable. Extending this to all
    arrows in $\pspath(r,s)$ we find that this latter description also
    succinctly summarises the action of $\pspath(\face_k)$ and
    immediately provides the following case-wise description of the
    related simplicial functor $\hcpath(\face_k)$ on a homset
    $\hcpath(r,s)$:
    \begin{itemize}
    \item $s<k$ wherein it is simply the identity on
      $\hcpath(r,s)$,
    \item $k\leq r$ when it is the canonical isomorphism between
      $\hcpath(r,s)$ and $\hcpath(r+1,s+1)$,
    \item $r<k\leq s$ in which case it is the simplicial map which
      carries a simplex $(\alpha_s,...,\alpha_{r+1})$ of
      $\hcpath(r,s)$ to the simplex
      $(\alpha_s,...,\alpha_k,1,\alpha_{k-1},...,\alpha_{r+1})$ in
      $\hcpath(r,s+1)$, or in other words we can picture it as the
      restriction:
      \begin{equation*}
        \xymatrix@R=2em@C=14em{
          {\hcpath(r,s)}
          \ar[r]^-{\hcpath(\face_k)}\ar@{u(->}[d]_{\subseteq_r} &
          {\hcpath(r,s+1)}
          \ar@{u(->}[d]^{\subseteq_r} \\
          {\Delta[1]^{(s-r)}}
          \ar[r]_-{\Delta[1]^{(s-k)}\times\Delta(\vertex_1)\times
          \Delta[1]^{(k-r-1)}} & {\Delta[1]^{(s-r+1)}}}
      \end{equation*}
    \end{itemize}
  \item {\boldmath $\arrow
      \hcpath(\degen_k):\hcpath[n+1]->\hcpath[n].$} Again this acts on
    objects $r\in[n+1]$ by applying $\degen_k$ and we may again
    determine its action on homsets by considering how the
    corresponding functor $\pspath(\face_k)$ acts on an indecomposable
    $\indec<r,s>$:
    \begin{equation*}
      \pspath(\degen_k)\indec<r,s> =
      \begin{cases}
        \indec<r,s> & \text{if $s< k+1$,} \\
        \indec<r-1,s-1> & \text{if $k+1\leq r$,} \\
        \indec<r,s-1> & \text{if $r< k+1 \leq s$.}
      \end{cases}
    \end{equation*}
    This time we cannot simply interpret this as the mere removal of a
    symbol at index $k+1$, since if $k+1=s$ this would result in the
    removal of a terminating $0$ to give a tuple which is not an
    element of $\pspath(r,s-1)$. A moment's reflection, however,
    reveals that we may instead summarise the action of
    $\pspath(\degen_k)$ on a general arrow $(a_s,a_{s-1},...,a_{r+1})$
    of $\pspath(r,s)$ by saying that it drops the rightmost ordinate
    $a_{r+1}$ if $r=k$ and otherwise replaces the ordinates $a_k$ and
    $a_{k+1}$ by a single entry obtained by taking the minimum
    $\min(a_k,a_{k+1})$ of these two values if $r<k<s$. We may now
    apply this observation to obtain a case-wise description of
    $\hcpath(\degen_k)$ on a homset $\hcpath(r,s)$:
    \begin{itemize}
    \item $s\leq k$ wherein it is simply the identity on $\hcpath(r,s)$,
    \item $k<r$ when it is the canonical isomorphism between
      $\hcpath(r,s)$ and $\hcpath(r-1,s-1)$,
    \item $k=r$ in which case it is the map which drops the initial
      ordinate of each simplex $(\alpha_s,...,\alpha_{r+1})$ in
      $\hcpath(r,s)$, or in other words we can picture it as the
      restriction:
      \begin{equation*}
        \xymatrix@R=2em@C=10em{
          {\hcpath(r,s)}
          \ar[r]^-{\hcpath(\degen_k)}\ar@{u(->}[d]_{\subseteq_r} &
          {\hcpath(r,s-1)}
          \ar@{u(->}[d]^{\subseteq_r} \\
          {\Delta[1]^{(s-r)}}
          \ar[r]_-{\Delta[1]^{(s-r-1)}\times\Delta(\degen^0_0)}
          & {\Delta[1]^{(s-r-1)}}}
      \end{equation*}
    \item $r<k<s$ in which case it is the map which replaces ordinates
      $\alpha_k$ and $\alpha_{k+1}$ of each simplex
      $(\alpha_s,...,\alpha_{r+1})$ in $\hcpath(r,s)$ by their
      pointwise minimum $\min(\alpha_k,\alpha_{k+1})$, or in other
      words we can picture it as the restriction:
      \begin{equation*}
        \xymatrix@R=2em@C=14em{
          {\hcpath(r,s)}
          \ar[r]^-{\hcpath(\degen_k)}\ar@{u(->}[d]_{\subseteq_r} &
          {\hcpath(r,s-1)}
          \ar@{u(->}[d]^{\subseteq_r} \\
          {\Delta[1]^{(s-r)}}
          \ar[r]_-{\Delta[1]^{(s-k-1)}\times\min\times\Delta[1]^{(k-r-1)}}
          & {\Delta[1]^{(s-r-1)}}}
      \end{equation*}
    \end{itemize}
  \end{itemize}
  Consequently, in each of these cases we may apply the argument
  outlined in the first paragraph, thereby showing that
  $\hcpath(\face_k)$ and $\hcpath(\degen_k)$ lift to
  $\Slgray$-enriched functors $\arrow\hclgray(\face_k):
  \hclgray[n-1] ->\hclgray[n].$ and $\arrow
  \hclgray(\degen_k):\hclgray[n+1]->\hclgray[n].$,
  and thus completing the demonstration that the functor $\hcpath$
  lifts as promised.
\end{obs}

\begin{notation}
  Let $\sint(r,s)$ denote $\{i\in\mathbb{N}\mid r<i<s\}$ the strict
  interval of integers between $r$ and $s$ and let $\hint(r,s]$ denote
  $\{i\in\mathbb{N}\mid r<i\leq s\}$ the corresponding half strict
  interval. Our primary use for these intervals is that they index the
  ordinates of simplices in an iterated product $\Delta[1]^n$ or a
  homset $\hclgray(r,s)$.

  It is clear that the $m$-simplices of $\Delta[1]^{(s-r)}$ correspond
  to arbitrary functions $\arrow w:\hint(r,s]->\cube{m}.$, simply
  because the operators of definition~\ref{simp.of.1.a} enumerate all of
  the $m$-simplices of $\Delta[1]$ and are indexed by $\cube{m}$. The
  simplex $\vec\rho_w$ associated with such a $w$ is simply the
  $(s-r)$-tuple $(\rho^m_{w(s)}, \rho^m_{w(s-1)}, ...,
  \rho^m_{w(r+1)})$.  We will, if necessary, implicitly extend a
  function $\arrow w:\sint(r,s)->\cube{m}.$ to a function with domain
  $\hint(r,s]$ by letting $w(s)={-}$, under which convention the
  corresponding $m$-simplex $\vec\rho_w$ becomes an element of the
  homset $\hclgray(r,s)$. Furthermore, we say that such a function is:
  \begin{enumerate}[label=(\roman*)]
  \item a {\em partial bijection\/} if for each integer $i\in\cube{m}$
    there exists a unique integer $j$ in its domain such that $w(j)=i$.
  \item {\em order reversing\/} if whenever $i,j$ are integers in the
    domain of $w$ with $i\leq j$ and for which $f(i)$ and $f(j)$ are
    integers in $\cube{m}$ then $f(j)\leq f(i)$.
  \end{enumerate}
  Notice that these definitions impose no conditions at those integers
  in the domain of $w$ which actually map to one of the points $-$ or
  $+$ in its codomain.
\end{notation}

\begin{obs}[simplicial cubes]\label{simp.cube.obs}
  In order to understand the homsets $\hcpath(r,s)$ more thoroughly,
  we should rehearse some basic facts about the structure of the {\em
    simplicial cube\/} $\Delta[1]^n$ and apply them to analysing our
  {\em stratified cube\/} $\Delta[1]^{\laxgray n}$.  Firstly, it is
  clear that the $m$-simplex $\vec\rho_w\in\Delta[1]^n$ associated
  with a function $\arrow w:\hint(0,n]->\cube{m}.$ is degenerate iff
  there is some integer $k\in[m-1]$ such that
  $\rho^r_{w(i)}(k)=\rho^r_{w(i)}(k+1)$ for all $i=1,2,...,n$ and
  this happens iff $\rho^r_{k+1}$ is not numbered amongst its
  ordinates, since this is the only such operator which maps $k$ and
  $k+1$ to different elements of $[1]$. In other words, the simplex
  corresponding to $w$ is degenerate iff $w$ is not surjective onto
  the integers (non-points) of $\cube{m}$.  This immediately implies
  that all of the simplices of $\Delta[1]^n$ of dimension
  greater than $n$ are degenerate and that its non-degenerate
  $n$-simplices correspond to partial bijections $\arrow
  w:\hint(0,n]->\cube{n}.$.

  Extending this analysis to study the stratification of
  $\Delta[1]^{\laxgray n}$, we may show that if $\arrow w:\hint(0,n]
  ->\cube{m}.$ is a partial bijection then the corresponding
  $m$-simplex $\vec\rho_w$ is non-thin if and only if $w$ is order
  reversing. First observe that we only need demonstrate that the
  result holds when $n=m$, since if $m<n$ then $\vec\rho_w$ is a
  special simplex in some $m$-dimensional cubical face of
  $\Delta[1]^{\laxgray n}$ isomorphic to $\Delta^{\laxgray m}$. To
  prove the ``only if'' part, we will assume, for a contradiction,
  that there exists some pair of integers $i<j$ in $\hint(0,n]$ for
  which $w(i)<w(j)$. Of course, we know that $\Delta[1]^{\laxgray n}$
  is isomorphic to $\Delta[1]^{\laxgray
    (n-i)}\laxgray\Delta[1]^{\laxgray i}$ and if we express
  $\vec\rho_w$ as a simplex in the latter set it becomes a pair with
  first component $(\rho^n_{w(n)},...,\rho^n_{w(i+1)})$ and second
  component $(\rho^n_{w(i)},...,\rho^n_{w(1)})$. Since $w$ is a
  partial bijection we know that $\rho^n_{w(i)}$ does not occur as an
  ordinate of the first of these components, (since it occurs in the
  second) and thus that it is degenerate at $w(i)-1$ by the argument
  of the last paragraph. Dually we see that the second component is
  degenerate at $w(j)-1$. Now we may apply lemma~129
  of~\cite{Verity:2006:Complicial} to show that our pair
  $((\rho^n_{w(n)},...,\rho^n_{w(i+1)}),(\rho^n_{w(i)},
  ...,\rho^n_{w(1)}))$ is thin in $\Delta[1]^{\laxgray
    (n-i)}\laxgray\Delta[1]^{\laxgray i}$ and thus that $\vec{\rho}_w$
  is thin in $\Delta[1]^{\laxgray n}$, thereby providing the desired
  contradiction.

  For the implication in the ``if'' direction we define an order
  preserving function $\arrow \bar{c}^n:[1]^n->[n].$ by letting
  \begin{equation*}
    \bar{c}^n(a_n,a_{n-1},...,a_1) \defeq \min( \{n\}\cup\{n-i\mid
    i\in\hint(0,n] \wedge a_i=0\}) 
  \end{equation*}
  from which we construct a simplicial map $\arrow
  c^n:\Delta[1]^{n}->\Delta[n].$ by applying the categorical nerve
  construction $\arrow \nerv:\Cat->\Simp.$. It is the case that this is
  the underlying simplicial map of a stratified map $\arrow
  c^n:\Delta[1]^{\laxgray n}->\Delta[n].$, which fact we prove via a
  simple induction which re-expresses $\bar{c}^n$ as a composite
  \begin{equation*}
    \xymatrix@R=1ex@C=8em{
      {[1]^n}\ar[r]^{[1]\times \bar{c}^{n-1}} &
      {[1]\times[n-1]}\ar[r]^{\bar{d}^n} & {[n]}
    }
  \end{equation*}
  wherein the order preserving map $\bar{d}^n$  is given by:
  \begin{equation*}
    \bar{d}^n(i,j) =
    \begin{cases}
      0 & \text{if $i=0$,} \\
      j+1 & \text{if $i=1$.}
    \end{cases}
  \end{equation*}
  Applying the nerve construction to this composite we obtain a pair
  of simplicial maps which compose to give $c^n$. Considering these
  factors in turn, we know that $c^{n-1}$ underlies a suitable
  stratified map, by the induction hypothesis, so the first of these
  $\Delta[1]\times c^{n-1}$ is the underlying simplicial map of the
  stratified map $\arrow \Delta[1]\laxgray c^{n-1}:\Delta[1]^{\laxgray
    n}-> \Delta[1]\laxgray \Delta[n-1].$. To show that the second
  factor also underlies a suitable stratified map, we may apply
  corollary~130 of loc.\ cit.\ to find that a
  simplex $(\alpha,\beta)\in\Delta[1]\laxgray\Delta[n-1]$ is thin if
  and only if there exists a $k\leq l$ such that
  $\alpha(k)=\alpha(k+1)$ and $\beta(l)=\beta(l+1)$. If that thin
  simplex is also non-degenerate then it follows that we have $0\leq
  k<l$ and $\alpha(l)<\alpha(l+1)$, so $\alpha(l-1)=\alpha(l)=0$ and
  it follows that $\bar{d}^n$ maps both of the vertices
  $(\alpha(l-1),\beta(l-1))$ and $(\alpha(l),\beta(l))$ to the vertex
  $0$ of $\Delta[n]$. Consequently, $d^n$ maps our thin simplex
  $(\alpha,\beta)$ to a simplex in $\Delta[n]$ which is degenerate at
  $(l-1)$, and is thus thin, from which it follows that $d^n$
  also underlies a stratified map $\arrow d^n:
  \Delta[1]\laxgray\Delta[n-1]->\Delta[n].$. So we
  find that the composite of these factors is a stratified map $\arrow
  c^n:\Delta[1]^{\laxgray n}->\Delta[n].$ and we may now prove that
  the $n$-simplex $\vec\rho_w$ is non-thin in
  $\Delta[1]^{\laxgray n}$ simply by observing that $c^n$ maps it to
  the unique non-degenerate $n$-simplex $\id_{[n]}$ in $\Delta[n]$
  which is non-thin in there.
\end{obs}

\begin{defn}
  We shall say that a function $\arrow w:\sint(r,s)->\cube{m}.$ is
  {\em special\/} if it is partially bijective and order
  reversing. Consulting observation~\ref{simp.cube.obs} we see that
  these conditions ensure that the corresponding special simplex
  $\vec{\rho}_w$ is a non-thin $m$-simplex of
  $\hclgray(r,s)$. Furthermore, this special simplex is an
  indecomposable arrow of $\hclgray(r,s)$ if and only if $w(i)\neq{-}$
  for all integers $i\in\sint(r,s)$, in which case we shall say that
  $w$ itself is indecomposable. We reserve the notation $\arrow
  s_n:\sint(0,n+1)->\cube{n}.$ to denote the unique such special
  function, given by $s_n(i)=n-i+1$.
\end{defn}

\begin{defn}[the nerve of a Gray-category]
  Applying Kan's construction to the functor $\arrow\hclgray:
  \Delta->\ECat{\Slgray}.$ we obtain an adjoint pair:
  \begin{equation*}
    \adjdisplay \ECat{\Slgray}->\Simp, \ladj-| \nerv.
  \end{equation*}
  To stratify the nerve of a $\Slgray$-category $\lcat{E}$ above
  dimension $1$ we specify that an $n$-simplex $\arrow
  f:\hclgray[n]->\lcat{E}.$ (with $n>1$) will be taken to be thin in
  the nerve $\nerv(\lcat{E})$ iff it maps the unique non-degenerate
  $(n-1)$-simplex of the homset $\hclgray(0,n)$ to a thin simplex in
  the stratified set $\lcat{E}(f(0),f(n))$. We can then specify that a
  $1$-simplex will be thin in $\nerv(\lcat{E})$ iff it has an equivalence
  inverse in the sense of theorem~\ref{almost.weakinner.cor}
  of~\cite{Verity:2005:WeakComp} with respect to the given
  stratification at dimension $2$.  It is now a matter of routine
  verification to check that this stratification is well defined and
  functorial.

  Now we know that the family of entire inclusions
  $\overinc\subseteq_e:X\laxgray Y->X\gray Y.$ make the identity on
  $\Strat$ into a monoidal functor from $(\Strat,\gray,\Delta[0])$ to
  $(\Strat,\laxgray,\Delta[0])$.  This allows us to regard every
  $\Strat\sgray$-category as a $\Slgray$-category, and it
  follows that we may apply the nerve construction above to any
  Gray-category. We shall show that the nerve of a Gray-category
  $\lcat{G}$ is always a weak complicial set, for which we only need
  check that $\nerv(\lcat{G})$ has fillers for inner horns and then apply
  theorem~\ref{almost.weakinner.cor} of~\cite{Verity:2005:WeakComp}.
\end{defn}

\begin{lemma}
  The nerve functor $\arrow \nerv:\ECat\Slgray->\Simp.$ is faithful.
\end{lemma}

\begin{proof}
  To prove this result we start by defining a locally ordered category
  $\Sigma[n]$ which has two objects $0$ and $1$ and homsets
  $\Sigma[n](0,1)\defeq[n]$ and
  $\Sigma[n](0,0)=\Sigma[n](1,1)\defeq[0]$ under the obvious (trivial)
  composition. Now we may define a local order preserving functor
  $\arrow \bar{f}^n:\pspath[n+1]->\Sigma[n].$ which maps the objects
  $0,1,...,n$ to $0$, the object $n+1$ to $1$ and acts on an
  arrow $\vec{a}$ in the homset $\pspath(r,s)$ as follows:
  \begin{equation*}
    \bar{f}^n(a_s,...,a_{r+1}) =
    \begin{cases}
      0 & \text{if $s\leq n$,} \\
      0 & \text{if $r,s=n+1$,} \\
      \min(\{n\}\cup\{n-i\mid i\in \hint(r,n]\wedge a_i=0\}) &
      \text{otherwise.}
    \end{cases}
  \end{equation*}
  To check that this is indeed functorial, simply observe that if
  $\vec{a}\in\pspath(s,n)$ and $\vec{b}\in\pspath(r,s)$ then their
  composite $\vec{a}\circ\vec{b}$ has a $0$ at index $s$ which ensures
  that the minima used to define
  $\bar{f}^n(\vec{a})=\bar{f}^n(\vec{a})\circ \bar{f}^n(\vec{b})$ and
  $\bar{f}^n(\vec{a}\circ\vec{b})$ coincide. Applying the categorical
  nerve construction $\arrow \nerv:\Cat->\Simp.$ to homsets we obtain a
  corresponding simplicial functor $\arrow
  f^n:\hcpath[n+1]->\Sigma\Delta[n].$, where the latter category has
  $\Delta[n]$ as its only non-trivial homset
  $\Sigma\Delta[n](0,1)$. However, by construction the action of $f^n$
  on the homset $\hcpath(0,n)$ is simply (isomorphic to) the
  simplicial map $\arrow c^n:\Delta[1]^n->\Delta[n].$ of
  observation~\ref{simp.cube.obs}, which we know to underlie a suitably
  stratified map thereby demonstrating that $f^n$ actually provides a
  $\Slgray$-functor $\arrow f^n:\hclgray[n+1]->\Sigma\Delta[n].$.

  Now observe that if $\lcat{E}$ is a $\Slgray$-category then (by
  Yoneda's lemma) a $\Slgray$-functor $\arrow
  g:\Sigma\Delta[n]->\lcat{E}.$ corresponds uniquely to an $n$-arrow
  in some homset of $\lcat{E}$. Furthermore if we compose this
  enriched functor with $\arrow f^n:\hclgray[n+1]->\Sigma\Delta[n].$
  we obtain an $(n+1)$-simplex $g\circ f^n$ of $\nerv(\lcat{E})$ from
  which we may regain the $n$-arrow that $g$ corresponds to by
  applying $g\circ f^n$ to the unique special $n$-simplex
  $\vec\rho_{s_n}$ in $\hclgray(0,n+1)\cong\Delta[1]^{\laxgray
    n}$. This fact is enough to show that $\nerv$ is faithful, since if
  $\arrow u,v:\lcat{C}->\lcat{D}.$ are two $\Slgray$-functors for
  which $\nerv(u)=\nerv(v)$ and we consider an arbitrary $n$-arrow $x$ of
  $\lcat{C}$ whose corresponding enriched functor we denote by
  $\arrow\yoneda{x}:\Sigma\Delta[n]->\lcat{C}.$ and for which we know
  that $u(x)=u(\yoneda{x}\circ f^n(\vec\psi_{s_n})) =
  \nerv(u)(\yoneda{x}\circ f^n)(\vec\psi_{s_n}) = \nerv(v)(\yoneda{x}\circ
  f^n)(\vec\psi_{s_n}) =v(\yoneda{x}\circ f^n(\vec\psi_{s_n})) =
  v(x)$. Consequently we have succeeded in demonstrating that $u$ and
  $v$ act identically on all arrows of $\lcat{C}$ and thus that they
  are identical as enriched functors.
\end{proof}

\section{Nerves of Gray Categories as Weak Complicial Sets}

\begin{obs}[complicial simplices in nerves of
  $\protect\Slgray$-categories]
  By definition, an inner complicial simplex $\overarr
  :\Delta^k[n]->\nerv(\lcat{E}).$ in the nerve of a $\Slgray$-category
  $\lcat{E}$ corresponds to an enriched functor $\arrow
  f:\hclgray[n]->\lcat{E}.$ satisfying an appropriate thinness
  property. To be precise, if $\arrow\alpha:[m]->[n].$ is a
  $k$-admissible operator then the composite $\arrow
  f\circ\hclgray(\alpha):\hclgray[m] ->\lcat{E}.$ must map the unique
  special $(m-1)$-simplex of $\hclgray(0,m)$ to a thin simplex in the
  homset $\lcat{E}(f(\alpha(0)), f(\alpha(m)))$.

  Arguing as we did in observation~\ref{hcpath.lift.obs},
  we may summarise the action of $\hclgray(\alpha)$, for any face
  operator $\arrow\alpha:[m]->[n].$, as a process of ``inserting the
  constant operator $1$ at indices in the set
  $[n]\setminus\im(\alpha)$''. More formally, if we define an
  associated indecomposable special function
  $\arrow\hat{\alpha}:\sint(\alpha(0), \alpha(m))-> \cube{m-1}.$ by
  \begin{equation*}
    \hat{\alpha}(i) =
    \begin{cases}
      j & \text{if $j\in\cube{m-1}$ and $\alpha(j)=i$,} \\
      + & \text{otherwise.}
    \end{cases}
  \end{equation*}
  then $\hclgray(\alpha)$ acts on the unique special $(m-1)$-simplex of
  $\hclgray(0,m)$ to map it to the corresponding special simplex
  $\vec{\rho}_{\hat\alpha}$ in $\hclgray(\alpha(0),\alpha(m))$.
  Observe also that any indecomposable special function $\arrow
  w:\sint(r,s)->\cube{m-1}.$ gives rise to a face operator $\arrow
  \check{w}:[m]->[n].$ given by
  \begin{equation*}
    \check{w}(j) =
    \begin{cases}
      r & \text{if $j=0$,} \\
      i & \text{if $w(i)=j$} \\
      s & \text{if $j=m$}
    \end{cases}
  \end{equation*}
  and that these two constructions are mutually inverse. Furthermore,
  the face operator $\check{w}$ is $k$-admissible if and only if
  $k\in\sint(r,s)$ and $w(i)$ is an integer (not a point) for all
  $i\in\{k-1,k,k+1\}\cap\sint(r,s)$, in which case we simply say that
  $w$ itself is $k$-admissible. Finally we see that an enriched
  functor $\arrow f:\hclgray[n]->\lcat{E}.$ corresponds to an inner
  $k$-complicial $n$-simplex in $\nerv(\lcat{E})$ if and only if it maps
  the simplex $\vec{\rho}_w$ to a thin simplex in $\lcat{E}$ for each
  $k$-admissible function $w$ as above. 
\end{obs}

\begin{obs}[complicial simplices in nerves of Gray-categories]
  \label{comp.simp.in.gray}
  Now if $\lcat{G}$ is actually a Gray-category then we find that any
  enriched functor $\arrow f:\hclgray[n]->\lcat{G}.$ which corresponds
  to a $k$-complicial simplex, as in the last observation, actually
  maps many of the decomposable simplices of $\hclgray[n]$ to thin
  simplices in the homsets of $\lcat{G}$. In particular if each
  component of the composite $\vec{\rho}_2\circ\vec{\rho}_1$ in
  $\hclgray$ maps to a thin arrow in $\lcat{G}$ then the pair
  $(f(\vec{\rho}_2),f(\vec{\rho}_1))$ is a thin simplex in
  $\lcat{G}(f(r),f(s))\gray\lcat{G}(f(s),f(t))$ which must thus map to
  a thin simplex $f(\vec{\rho}_2)\circ f(\vec{\rho}_1)=
  f(\vec{\rho}_2\circ\vec{\rho}_1)$ in $\lcat{G}(f(r),f(t))$. This is
  not generally the case for general $\Slgray$-categories, which
  prevents their nerves from being weak complicial sets even if their
  homsets happen to be such.

  Applying this to the result of the last observation, we find that
  any $k$-admissible simplex $\arrow f:\hclgray[n]->\lcat{G}.$ also
  maps the special simplex $\vec{\rho}_w$ associated with the function
  $\arrow w:\sint(r,s)->\cube{m-1}.$ to a thin simplex in a homset of
  $\lcat{G}$ if:
  \begin{enumerate}[label=\textbf{case (\roman*)},
    fullwidth, leftmargin=1em, itemsep=1ex]
  \item There are integers $i$,$j$ and $l$ in $\sint(r,s)$ with
    $i<l<j$ such that $w(l)={-}$ and both of $w(i)$ and $w(j)$ are
    integers (not points).
    \begin{proof}
      Observe that the condition $w(l)={-}$ implies that $\vec\rho_w$
      may be decomposed at index $l$ to express it as a composite
      $\vec\rho_2\circ\vec\rho_1$ where $\vec\rho_1\in\hclgray(r,l)$
      has the operator $\rho_{w(i)}$ as its $i\oth$ ordinate, which is
      non-constant since $w(i)$ is an integer, and similarly
      $\vec\rho_2\in\hclgray(l,s)$ has the non-constant $\rho_{w(j)}$
      as its $(j-l)\oth$ ordinate. However $w$ is partially bijective,
      so we know that each non-constant $\rho$ operator occurs exactly
      once as an ordinate of $\vec\rho_w$ and it follows that
      $\vec\rho_1$ does not have $\rho_{w(j)}$ as an ordinate and thus
      that it must be degenerate at $w(j)-1$. Dually $\vec\rho_2$ is
      degenerate at $w(i)-1$, so both of $f(\vec\rho_1)$ and
      $f(\vec\rho_2)$ are degenerate, and thus thin, simplices of
      their respective homsets and thus we find that $f(\vec\rho_w)$,
      their composite, is also thin in $\lcat{G}(r,s)$.
    \end{proof}
  \item We have $k\in\sint(r,s)$ and $w(k)$ is an integer and $w(i)$
    is an integer or the point $-$ (in other words is not the point
    $+$) for any $i\in\{k-1,k+1\}\cap\sint(r,s)$.
    \begin{proof}
      If $\vec\rho_w$ is not already indecomposable then we may
      decompose it as $\vec\rho_3\circ\vec\rho_2\circ\vec\rho_1$ in such
      a way that $\vec\rho_2$ is indecomposable with and ordinate
      indexed by $k$ and in which either of the other two factors
      could be identities. Now we have two possibilities, one of
      $\vec\rho_1$ or $\vec\rho_3$ has a non-constant ordinate, in
      which case we are reduced to the case presented above, or they
      are both degenerate and $\vec\rho_2$ is an indecomposable
      special simplex in some homset $\hclgray(r',s')$. Now it is
      clear that the doubly pointed morphism $\arrow
      w':\sint(r',s')->\cube{m-1}.$ corresponding to $\vec\rho_2$ is
      simply obtained by restricting $w$ to $\sint(r',s')$ and that
      the conditions given ensure that this $w'$ is $k$-admissible. So
      we know that $f(\vec\rho_2)$ is thin in its homset, since
      $\vec\rho_2$ is $k$-admissible, and that $f(\vec\rho_1)$ and
      $f(\vec\rho_3)$ are degenerate in their homsets, so it follows
      that their composite $f(\vec\rho_w)$ is thin in its homset as
      postulated.
    \end{proof}
  \end{enumerate}
\end{obs}

\begin{obs}
  Given this series of reductions, we are finally able to recast the
  task of proving that $\nerv(\lcat{G})$ is a weak complicial set as a
  demonstration that the members of a certain set of inclusions are
  complicial cofibrations. 

  Firstly our analysis of horns in homotopy coherent nerves, which
  culminated in observation~\ref{hcpath.of.graycat}, told us that in
  order to fill an inner horn in the nerve of homotopy coherent paths
  $\nerv_{hc}(\lcat{G})$ it is enough to demonstrate that each inclusion
  $\inc\subseteq_r:\hchorn^k[n](0,n)->\hcpath[n](0,n).$ is an anodyne
  extension. In the current case, it follows that all we need do is
  re-stratify $\hcpath[n](0,n)$ and $\hchorn^k[n](0,n)$ in line with
  the analysis of complicial simplices in $\nerv(\lcat{G})$ provided by
  observation~\ref{comp.simp.in.gray} and demonstrate that the
  corresponding regular inclusion is a complicial cofibration. Well in
  fact that isn't quite enough on its own, since this time the
  stratification of $\nerv(\lcat{G})$ is non-trivial at all dimensions so
  we will also need to ensure that the corresponding result holds for
  elementary inner thinness extensions. This however requires only a
  minor modification of the argument for horns.
\end{obs}

\begin{notation}\label{ckn.not}
  So to fix our notation we define the following stratified sets:
  \begin{itemize}
  \item $C^k_n$ ($n\geq 2$, $1 \leq k \leq n$) by starting with the
    stratified cube $\Delta[1]^{\laxgray n}$ and making thin those
    simplices $\vec\rho_w$ for which:
    \begin{enumerate}[label=(\roman*)]
    \item\label{ckn.crit.1} there are integers $i,j,l\in\hint(0,n]$
      with $i<l<j$ such that $w(l)=-$ and both of $w(i)$ and $w(j)$
      are integers,
    \item\label{ckn.crit.2} the value $w(k)$ is an integer and $w(i)$
      is an integer or the point $-$ for each
      $i\in\{k-1,k+1\}\cap\hint(0,n]$.
    \end{enumerate}
  \item $H^k_n$ is the regular subset of $C^k_n$ whose underlying
    simplicial set is the domain of the corner product depicted in
    display~(\ref{horn.corner.prod}).
  \item $\arrow w^n_i:\hint(0,n]->\cube{n-1}.$ ($1\leq i\leq n$) is
    the unique order reversing partial bijection for which
    $w_i(i)=+$. As ever we will habitually drop the superscript in
    this notation if it may be inferred from the context.
  \item $\dot{C}^k_n$ is obtained from $C^k_n$ by making the special
    $(n-1)$-simplices $\vec\rho_{w_i}$ thin for all
    $i\in\{k-1,k+1\}\cap\hint(0,n]$ and $\ddot{C}^k_n$ is obtained
    from that by making $\vec\rho_{w_k}$ thin as well.
  \end{itemize}
\end{notation}

\begin{obs}\label{proof.outline.obs}
  Of course, observation~\ref{comp.simp.in.gray} tells us that
  $C^k_{n-1}$ is simply the set obtained by taking $\hclgray(0,n)$ and
  making thin those simplices that would be mapped to thin simplices
  by any enriched functor $\arrow f:\hclgray[n]-> \lcat{G}.$ that
  corresponds to a $k$-complicial $n$-simplex in the nerve of the
  Gray-category $\lcat{G}$. Furthermore, $H^k_{n-1}$ is the stratified
  set with the corresponding property for $k$-complicial horns in
  $\nerv(\lcat{G})$. It follows, therefore, that we may demonstrate that
  $\nerv(\lcat{G})$ has fillers for such horns by showing that the
  inclusion $\overinc \subseteq_r:H^k_n->C^k_n.$ is a complicial
  cofibration for all $n\geq 2$ and $1\leq k\leq n$. Correspondingly,
  we may demonstrate the thinness extension property required of a
  weak inner complicial set by showing that each inclusion
  $\overinc\subseteq_e:\dot{C}^k_n->\ddot{C}^k_n.$ is a complicial
  cofibration.
\end{obs}

\begin{lemma}
  If $k>1$ then the Gray pre-tensor $C^k_{n-1}\pretens \Delta[1]$ is an
  entire subset of $C^{k+1}_n$ and the inclusion
  $\overinc\subseteq_e: C^k_{n-1}\pretens \Delta[1]->C^{k+1}_n.$ is an
  anodyne extension. Dually if $k<n-1$ then the corresponding result
  relates $\Delta[1]\pretens C^k_{n-1}$ and $C^k_n$.
\end{lemma}

\begin{proof}
  By definition we know that
  $\Delta[1]^{\laxgray(n-1)}\pretens\Delta[1]$ is an entire subset of
  $C^k_{n-1}\pretens\Delta[1]$ and, consulting
  definition~\ref{pretens.defn} of~\cite{Verity:2005:WeakComp}, we see
  that the first of the criteria listed there to define the
  stratifications of these pre-tensors depends only on underlying
  simplicial sets which in this case are identical. So applying the
  second of those criteria we see that an $m$-simplex
  $\vec\rho_w=(\vec\rho,\rho_{w(1)})$ is thin in the latter set and
  non-thin in the former if and only if either
  \begin{itemize}
  \item $w(1)$ is equal to $-$ or $+$ and $\vec\rho$ is special and
    satisfies one of the criteria for thinness given in the definition
    of $C^k_{n-1}$, or
  \item $w(1)$ is equal to $m$ (thus ensuring that
    $\rho_{w(1)}=\partproj^{m-1,1}_2$), the $(m-1)$-simplex
    $\vec\rho\cdot\partinj^{m-1,1}_1$ is special and satisfies
    one of the criteria given in the definition of $C^k_{n-1}$ and
    $(\vec\rho\cdot\partinj^{m-1,1}_1)\cdot\partproj^{m-1,1}_1 =
    \vec\rho$
  \end{itemize}
  with all other possibilities for the application of that thinness
  criterion turning out to be degenerate simplices. 

  In the first case, if $\arrow w':\sint(0,n)->\cube{m}.$ is the
  unique special function with $\vec\rho=\vec\rho_{w'}$ then it is
  clear that $w(i+1)=w'(i)$ and so $w$ is also special since
  $w(1)=-/+$.  In the second case, if we now take $w'$ to be the
  special function corresponding to the simplex $\vec\rho
  \cdot \partinj^{m-1,1}_1$ we may replace the other equality given
  there to $\vec\rho_{w'} \cdot \partproj^{m-1,1}_1= \vec\rho$ from
  which we may again infer easily that $w(i+1)=w'(i)$ since
  $\rho^{m-1}_i\circ \partproj^{m-1,1}_1 = \rho^m_i$ for all
  $i\in\sint(0,m-1)$. So we see again that $w$ is special since $w'$
  maps to integers less than $w(1)=m$. Now it is a routine matter, in
  either case, to check that $w'$ satisfies one of the thinness
  conditions on $C^k_{n-1}$ given in observation~\ref{ckn.not} if and
  only if $w$ satisfies the corresponding thinness condition on
  $C^{k+1}_n$.  Consequently, since we know that $\Delta[1]^{\laxgray
    (n-1)}\pretens \Delta[1]$ is also an entire subset of $C^{k+1}_n$
  we have certainly shown that $C^k_{n-1}\pretens\Delta[1]$ is an
  entire subset of $C^{k+1}_n$. Indeed what we have actually shown is
  that $C^{k+1}_n=(C^k_{n-1}\pretens\Delta[1])\cup\Delta[1]^{\laxgray
    n}$ and thus that we have a pushout square:
  \begin{equation*}
    \xymatrix@R=2em@C=4em{
      {\Delta[1]^{\laxgray(n-1)}\pretens\Delta[1]}
      \ar@{u(->}[r]^-{\subseteq_e}
      \ar@{u(->}[d]_{\subseteq_e} & 
      {\Delta[1]^{\laxgray n}}\ar@{u(->}[d]^{\subseteq_e} \\
      {C^k_{n-1}\pretens\Delta[1]}\ar@{u(->}[r]_-{\subseteq_e} &
      {C^{k+1}_n}\poexcursion
    }
  \end{equation*}
  Now observation~\ref{pretens.gen.obs}
  of~\cite{Verity:2005:WeakComp}, which recollects lemma 139
  of~\cite{Verity:2006:Complicial}, tells us that the upper horizontal
  here is an anodyne extension so it follows that its lower horizontal
  is also an anodyne extension, as a pushout of such, as required.

  Notice that the condition that $k>1$ is vital to the argument of the
  last paragraph, since it ensures that the ``right hand flank'' of
  the pivotal ordinate $k$ is protected from the insertion of a $+$
  which might disrupt condition~\ref{ckn.crit.2} of
  observation~\ref{ckn.not}. In the case $k=1$ this result does not
  hold since, for instance, $(\rho_1,\rho_2,...,\rho_{n-1})$ is thin
  in $C^1_{n-1}$ but $(\rho_1,...,\rho_{n-1},+)$ is not thin in
  $C^2_n$.
\end{proof}

\begin{cor}\label{n.inner.horn.1}
  If each of the inclusions $\overinc\subseteq_r:H^1_2->C^1_2.$,
  $\overinc\subseteq_r:H^2_2->C^2_2.$ and
  $\overinc\subseteq_r:H^2_3->C^2_3.$ are complicial cofibrations then
  so are the inclusions $\overinc\subseteq_r: H^k_n->C^k_n.$
  and $\overinc\subseteq_e:\dot{C}^k_n->\ddot{C}^k_n.$ for
  $n\geq 2$ and $0<k\leq n$.
\end{cor}

\begin{proof}
  For $n\geq 3$ and $1<k<n$ consider the following commutative square
  \begin{equation*}
    \xymatrix@R=2em@C=6em{
      {(H^k_{n-1}\pretens \Delta[1])\cup(C^k_{n-1}\pretens\boundary\Delta[1])}
      \ar@{u(->}[r]^-{\subseteq_s}\ar@{u(->}[d]_{\subseteq_e} &
      {C^k_{n-1}\pretens\Delta[1]}
      \ar@{u(->}[d]^{\subseteq_e} \\
      {H^{k+1}_n}\ar@{u(->}[r]_-{\subseteq_r} &
      {C^{k+1}_n}
    }
  \end{equation*}
  in which the right hand vertical is an anodyne extension by the last
  lemma as, quite clearly, is the left hand vertical by a minor
  modification of the proof of that same lemma. So these verticals are
  both complicial cofibrations and we may apply the 2-of-3 property to
  show that the upper horizontal here is a complicial cofibration if
  and only if its lower horizontal is such. Observe, however, that the
  upper horizontal is actually the corner pre-tensor
  \begin{equation*}
    (\overinc\subseteq_r:H^k_{n-1}->C^k_{n-1}.)\cpretens
    (\overinc\subseteq_r:\boundary\Delta[1]->\Delta[1].)
  \end{equation*}
  which we know is a complicial cofibration whenever
  $\overinc\subseteq_r:H^k_{n-1}->C^k_{n-1}.$ is such, by
  theorem~\ref{ccfib.ctens.thm}. Applying this result inductively we
  have shown that if $k>1$ then $\overinc
  \subseteq_r:H^{k+1}_n->C^{k+1}_n.$ is a complicial cofibration
  whenever $\overinc\subseteq_r:H^2_{n-k+1}->C^2_{n-k+1}.$ is
  such. Furthermore, arguing dually, we see that if $k<n-1$ then
  $\overinc\subseteq_r:H^k_n->C^k_n.$ is a complicial cofibration
  whenever $\overinc\subseteq_r:H^k_{k+1}->C^k_{k+1}.$ is such.  So we
  may
  \begin{itemize}
  \item use the first of these reductions to show that the inclusion
    $\overinc\subseteq_r:H^n_n->C^n_n.$ is a complicial cofibration
    under the assumption that $\overinc\subseteq_r:H^2_2->C^2_2.$ is
    such, 
  \item use its dual to show that the inclusion
    $\overinc\subseteq_r:H^1_n->C^1_n.$ is a complicial cofibration
    under the assumption that $\overinc\subseteq_r:H^1_2->C^1_2.$ is
    such, and
  \item use each reduction in turn to show that the inclusion
    $\overinc\subseteq_r:H^k_n->C^k_n.$ with $1<k<n$ is a complicial
    cofibration assuming that $\overinc\subseteq_r:H^1_3->C^1_3.$ is
    such,
  \end{itemize}
  thereby demonstrating the desired result for all such
  inclusions. 

  To demonstrate the corresponding result for the inclusion
  $\overinc\subseteq_e:\dot{C}^k_n->\ddot{C}^k_n.$ start by observing
  that the $(n-1)$-simplex $\vec\rho_w\in C^k_n$ corresponding to some
  function $\arrow w:\hint(0,n]->\cube{n-1}.$ is in $H^k_n$ if and
  only if there is an $i\in\hint(0,n]$ for which $w(i)=-$ or ($i\neq
  k$ and $w(i)=+$). Consequently, by observation~\ref{simp.cube.obs},
  it follows that if such an $(n-1)$-simplex is non-degenerate then
  $w$ must be a partial bijection and if it is non-thin then $w$ must
  also be order reversing and thus special. Consulting the definition
  of $C^k_n$ we see that it only non-thin special $(n-1)$-simplices
  are the $\vec\rho_{w_i}$ for each $i\in\{k-1,k,k+1\}\cap\sint(1,n)$
  and that $\vec\rho_{w_k}$ is not a simplex of $H^k_n$. It follows,
  therefore, that $\dot{C}^k_n=C^k_n\cup\triv_{n-2}(H^k_n)$ and we may
  thus apply corollary~\ref{ccof.th.cor}
  of~\cite{Verity:2005:WeakComp} to the complicial cofibration
  $\overinc\subseteq_r:H^k_n->C^k_n.$ of the last paragraph to
  demonstrate that $\overinc\subseteq_e:\dot{C}^k_n->
  \triv_{n-2}(C^k_n).$ is a complicial cofibration. Now it is the case
  that if the composite of two {\em entire\/} inclusions is a
  complicial cofibration then so is each inclusion separately, which
  fact we leave to the reader to demonstrate using the
  characterisation given in corollary~\ref{ccof.char.b}
  of loc.\ cit. Applying this to the inclusion of
  the sentence before last, which factors as a composite of the entire
  inclusions $\overinc\subseteq_e:\dot{C}^k_n-> \ddot{C}^k_n.$ and
  $\overinc\subseteq_e:\ddot{C}^k_n-> \triv_{n-2}(C^k_n).$, we obtain
  the desired result.
\end{proof}

\begin{lemma}\label{n.inner.horn.2}
  The inclusions $\overinc\subseteq_r:H^1_2->C^1_2.$,
  $\overinc\subseteq_r:H^2_2->C^2_2.$ and
  $\overinc\subseteq_r:H^2_3->C^2_3.$ are indeed complicial
  cofibrations.
\end{lemma}

\begin{proof}
  The easiest way to visualise the calculations involved here is
  diagrammatically. The first two inclusions are duals of each other
  and the first may be pictured as follows
  \begin{equation*}
    \xymatrix@=2em{
      {(0,1)}\ar@{..>}[r] & {(1,1)} \\
      {(0,0)}\ar[u]|{\objectstyle\sim}\ar@{..>}[ur]\ar[r] & 
      {(0,1)}\ar[u]
    }
  \end{equation*}
  wherein the solid arrows represent the $1$-simplices of $C^1_2$
  which are in $H^1_2$ and the dotted ones represent $1$-simplices
  which will be obtained by filling some $2$-simplex. The left hand
  vertical $1$-simplex is labelled with a $\sim$ symbol because it is
  thin in $C^1_2$ by the virtue of satisfying the first condition in
  the definition of the stratification of $C^2_3$ given in
  notation~\ref{ckn.not}. We may decompose the regular subset
  inclusion $\overinc\subseteq_r:H^1_2->C^1_2.$ by defining the
  intermediate regular subset
  \begin{equation*}
    U \defeq H^1_2\cup \{(0,0)<(0,1)<(1,1)\}^*
  \end{equation*}
  where the notation $\{{-}\}^*$ denotes the regular subset generated
  by the given set of simplices, and show that:
  \begin{itemize}
  \item The inclusion $\overinc\subseteq_r:H^1_2->U.$ is anodyne since
    it is a pushout of the horn extension
    $\overinc\subseteq_r:\Lambda^1[2]->\Delta^1[2].$ along the evident
    stratified map mapping its domain to the $1$-complicial horn on
    the vertices $(0,0)$, $(0,1)$ and $(1,1)$.
  \item The inclusion $\overinc\subseteq_r:U->C^1_2.$ is anodyne since
    it is a pushout of the horn extension
    $\overinc\subseteq_r:\Lambda^0[2]->\Delta^1[2].$ along the evident
    stratified map mapping its domain to the $1$-complicial horn on
    the vertices $(0,0)$, $(1,0)$ and $(1,1)$.
  \end{itemize}
  It follows therefore that their composite
  $\overinc\subseteq_r:H^1_2->C^1_2.$ is an anodyne extension, and
  thus is a complicial cofibration as required.

  The proof for the inclusion $\overinc\subseteq_r:H^2_3->C^2_3.$ is
  only slightly more involved, again we assist the reader with a
  diagram
  \begin{equation}
    \label{cb.fill.2.d}
    \xymatrix@=2em{
      & {(0,1,1)}\ar[r] & {(1,1,1)} \\
      {(0,0,1)} \ar[ur]\ar@{}[urr]|{\objectstyle\vec\rho_{w_3}}\ar[r] &
      {(1,0,1)} \ar[ur] & {(1,1,0)}\ar[u] \\
      {(0,0,0)}\ar[u]\ar[r]\ar@{}[ur]|{\objectstyle\sim} &
      {(1,0,0)}\ar[u]\ar@{}[uur]|{\objectstyle\vec\rho_{w_1}}\ar[ur] & {}
    }\mkern40mu
    \xymatrix@=2em{
      & {(0,1,1)}\ar[r] & {(1,1,1)} \\
      {(0,0,1)} \ar[ur] &
      {(0,1,0)} \ar[u]\ar[r]\ar@{}[ur]|{\objectstyle?(\vec\rho_{w_2})} & {(1,1,0)}\ar[u] \\
      {(0,0,0)}\ar[u]\ar[r]\ar[ur]|{\objectstyle\sim} 
      \ar@{}[uur]|{\objectstyle\sim}\ar@{}[rru]|{\objectstyle\sim}&
      {(1,0,0)}\ar[ur] & {}
    }
  \end{equation}
  depicting the set $H^2_3$, Here the labels in the square faces and
  on one of the arrows provide information about the corresponding
  special $2$-simplices. Those labelled with $\sim$ symbols are thin
  in $H^2_3$, with the one in the left hand diagram satisfying the
  first condition in the definition of the stratification of $C^2_3$
  given in notation~\ref{ckn.not} and those in the right hand diagram
  all satisfying the second condition given there. The remaining
  labels simply identify the corresponding special simplex, with the
  one in the upper-right most square being prefixed with a $?$ since
  that is the square $2$-face which is not present in the set $H^2_3$.

  This time we actually need to consider the entire superset
  $\hat{C}^2_3$ of $C^2_3$ constructed by making the $2$-simplex
  $(0,0,0)<(0,1,0)<(1,1,1)$ thin and construct a decomposition of the
  inclusion $\overinc\subseteq_r:H^2_3->\hat{C}^2_3.$ consisting of an
  increasing tower of regular subsets of $\hat{C}^2_3$ given by:
  \begin{align*}
    V_1 & \defeq H^2_3\cup\{(0,0,0)<(0,1,1)<(1,1,1)\}^* \\
    V_2 & \defeq V_1\cup\{(0,0,0)<(0,0,1)<(0,1,1)<(1,1,1)\}^* \\
    V_3 & \defeq V_2\cup\{(0,0,0)<(0,0,1)<(1,0,1)<(1,1,1)\}^* \\
    V_4 & \defeq V_3\cup\{(0,0,0)<(1,0,0)<(1,0,1)<(1,1,1)\}^* \\
    V_5 & \defeq V_4\cup\{(0,0,0)<(1,0,0)<(1,1,0)<(1,1,1)\}^* \\
    V_6 & \defeq V_5\cup\{(0,1,0)<(0,1,1)<(1,1,1)\}^* \\
    V_7 & \defeq V_6\cup\{(0,0,0)<(0,1,0)<(0,1,1)<(1,1,1)\}^* \\
    \hat{C}^2_3 & = V_7\cup\{(0,0,0)<(0,1,0)<(1,1,0)<(1,1,1)\}^*
  \end{align*}
  Now we find that
  \begin{itemize}
  \item the inclusions $\overinc\subseteq_r:H^2_3->V_1.$ and
    $\overinc\subseteq_r:V_5->V_6.$ may be constructed as evident
    pushouts of the $1$-complicial horn
    $\overinc\subseteq_r:\Lambda^1[2]->\Delta^1[2].$,
  \item the inclusions $\overinc\subseteq_r:V_1->V_2.$ and
    $\overinc\subseteq_r:V_6->V_7.$ may be constructed as evident
    pushouts of the thin $2$-complicial horn
    $\overinc\subseteq_r:\Lambda^2[3]'->\Delta^2[3]''.$,
  \item the inclusion $\overinc\subseteq_r:V_3->V_4.$ may be
    constructed as a pushout of the $2$-complicial horn
    $\overinc\subseteq_r:\Lambda^2[3]->\Delta^2[3].$,
  \item the inclusions $\overinc\subseteq_r:V_2->V_3.$,
    $\overinc\subseteq_r:V_4->V_5.$ and
    $\overinc\subseteq_r:V_5->V_6.$ may be constructed as pushouts of
    the thin $1$-complicial horn
    $\overinc\subseteq_r:\Lambda^1[3]->\Delta^1[3].$, and
  \item the inclusion $\overinc\subseteq_r:V_7->\hat{C}^2_3.$ may be
    constructed as a pushout of the $0$-complicial horn
    $\overinc\subseteq_r:\Lambda^0[3]->\Delta^0[3].$.
  \end{itemize}
  as the reader may readily verify, so it follows that their composite
  $\overinc\subseteq_r:H^2_3->\hat{C}^2_3.$ is an anodyne
  extension. It is also easily demonstrated that the entire inclusion
  $\overinc\subseteq_e:C^2_3->\hat{C}^2_3.$ is an anodyne extension,
  since it may be constructed as a pushout of the elementary thinness
  extension $\overinc\subseteq_e:\Delta^2[3]'->\Delta^2[3]''.$ along
  the evident map of its domain onto the simplex
  $(0,0,0)<(0,0,1)<(0,1,1)<(1,1,1)$. So these inclusions are both
  complicial cofibrations, as are all anodyne extensions, so we may
  apply the 2-of-3 property to show that the inclusion
  $\overinc\subseteq_r:H^2_3->C^2_3.$ is also a complicial cofibration
  as required.
\end{proof}

\begin{defn}
  An enriched functor $\arrow f:\lcat{G}->\lcat{H}.$ of
  Gray-categories is said to be a {\em local complicial fibration\/}
  if each stratified map $\arrow
  f:\lcat{G}(a,b)->\lcat{H}(f(a),f(b)).$ of homsets is a complicial
  fibration.
\end{defn}

\begin{thm}
  If $\lcat{G}$ is a Gray-category then its nerve $\nerv(\lcat{G})$ is a
  weak complicial set. Furthermore, if $\arrow f:\lcat{G}->\lcat{H}.$
  is a functor of Gray-categories which is a local complicial
  fibration then the corresponding stratified map $\arrow
  \nerv(f):\nerv(\lcat{G})->\nerv(\lcat{H}).$ is an {\em inner\/} complicial
  fibration.
\end{thm}

\begin{proof}
  This is now just a simple matter of following the proof outline
  given in observation~\ref{proof.outline.obs} using the results
  established in corollary~\ref{n.inner.horn.1} and
  lemma~\ref{n.inner.horn.2} to demonstrate that $\nerv(\lcat{G})$ is
  (almost) a weak inner complicial set whenever $\lcat{G}$ is a
  Gray-category. Now we may apply theorem~\ref{almost.weakinner.cor}
  of~\cite{Verity:2005:WeakComp} to show that $\nerv(\lcat{G})$ is
  actually a weak complicial set as postulated. The remainder is a
  triviality, since it is clear that the results of
  corollary~\ref{n.inner.horn.1} and lemma~\ref{n.inner.horn.2} also
  show that we may lift inner horns and thinness extensions along
  $\arrow \nerv(f):\nerv(\lcat{G})->\nerv(\lcat{H}).$ so long as its actions on
  homsets has the RLP with respect to the inclusions featured in those
  lemmata.
\end{proof}

\begin{obs}[the large universe of weak complicial sets]
  We have now succeeded in showing that the weak complicial sets
  themselves are the $0$-simplices of a large and richly structured
  weak complicial set. This is simply constructed by taking the nerve
  of the Gray-category $\Wcs$ of weak complicial sets discussed in
  example~\ref{wcs.gray.obs.a}. It is worth observing that we may
  easily show that this is {\bf not} $n$-trivial for any
  $n\in\mathbb{N}$.  
\end{obs}

\bibliographystyle{amsplain}
\bibliography{cattheory}

\end{document}